\documentclass[letterpaper, 10 pt, conference]{ieeeconf}  

\IEEEoverridecommandlockouts                              

\usepackage{graphicx} 
\usepackage{amsmath} 
\usepackage{amssymb}  
\usepackage{tikz}

\usepackage{url}
\usepackage[pdfstartview=FitV, urlcolor=blue]{hyperref}

\newcommand{\bx}{\boldsymbol{x}}
\newcommand{\by}{\boldsymbol{y}}
\newcommand{\bz}{\boldsymbol{z}}
\newcommand{\ba}{\boldsymbol{a}}
\newcommand{\br}{\boldsymbol{r}}
\newcommand{\bh}{\boldsymbol{h}}
\newcommand{\xT}{\bx_{\textsc{t}}}
\newcommand{\xS}{\bx_{\textsc{s}}}
\newcommand{\x}{\bx}

\newcommand{\z}{\bz}
\newcommand{\A}{\mathcal{A}}
\newcommand{\I}{\mathcal{I}}
\newcommand{\J}{\mathcal{J}}
\newcommand{\Z}{\mathcal{Z}}
\newcommand{\V}{\mathcal{V}}
\newcommand{\E}{\mathbb{E}}

\DeclareMathOperator*{\argmin}{arg\,min}

\title{\LARGE \bf
Time-Dependent Surveillance-Evasion Games
}

\author{Elliot Cartee$^{1,2}$, Lexiao Lai$^{3}$, Qianli Song$^{3}$, Alexander Vladimirsky$^{1}$
\thanks{*Much of this work was conducted in an REU program, partially supported by the NSF-RTG award (DMS-1645643). The 1st and 4th authors were also supported by the NSF ATD award (DMS-1738010).  The last author's work is also supported by the Simons Foundation Fellowship.}
\thanks{$^{1}$Department of Mathematics, Cornell University}%
\thanks{$^{2}$evc34@cornell.edu}%
\thanks{$^{3}$University of Hong Kong}%
}

\begin{document}

\maketitle
\thispagestyle{empty}
\pagestyle{empty}

\begin{abstract}

Surveillance-Evasion (SE) games form an important class of adversarial trajectory-planning problems. 
We consider time-dependent SE games, in which an Evader is trying to reach its target while minimizing the cumulative exposure to a moving enemy Observer.  That Observer is simultaneously aiming to maximize the same exposure by choosing how often to use each of its predefined patrol trajectories.  Following the framework introduced in \cite{gilles2018surveillance}, we develop efficient algorithms for finding Nash Equilibrium policies for both players by blending techniques from semi-infinite game theory, convex optimization, and multi-objective dynamic programming on continuous planning spaces.
We illustrate our method on several examples with Observers using omnidirectional and angle-restricted sensors on a domain with occluding obstacles.    

\end{abstract}

\section{INTRODUCTION}

While optimal trajectory-planning is a common task in robotics, the notion of ``optimality'' can be based on many different criteria, such as minimizing time taken to reach the destination, energy required to traverse the path, or minimizing exposure to some sort of threat.
In this paper, we focus on the latter in the context of Surveillance-Evasion (SE) games.

We model this adversarial trajectory-planning problem as a semi-infinite zero-sum game between two robotic players: an Observer (O) and an Evader (E) traveling through a domain with occluding obstacles.   
O chooses a probability distribution over a finite set of predefined surveillance plans.
E chooses a probability distribution over an infinite set of trajectories that bring it to a target location before a specified deadline. O's position and orientation define a time-dependent {\em pointwise observability} function.
E's goal is to minimize its expected {\em cumulative observability}, while O aims to maximize this same quantity.
Our model is continuous in time and space, relying on numerical methods for partial differential equations (PDEs) to plan E's trajectories via dynamic programming.  We also use techniques from multi-objective and convex optimization to find Nash equilibrium policies for both players.

Classical SE games were introduced in the 1970s and assumed that each player is fully aware of the opponent's state and can react to any changes in real time \cite{LewinBreakwell,LewinOlsder}.  In contrast, our version of the game is built on a different information structure, as each player makes their decisions based only on the probabilistic policy chosen by its opponent, rather than on that opponent's actual current state.  In this sense, we are looking for optimal {\em open loop} controls, with a built-in uncertainty of the opponent's position. 
This information structure arises in applications where logistical considerations force the players to commit to strategies in advance.

This version of SE games was recently considered in \cite{gilles2018surveillance} for a model that was rather simplified from the robotics point of view: it lacked any detailed kinematics for E, allowed a set of stationary locations (rather than general surveillance plans) for O, and used the pointwise observability based on omnidirectional sensors at O's location(s).  Here we use the same algorithms for finding Nash equilibrium policies, but extend the approach introduced in \cite{gilles2018surveillance} to more realistic settings.  While we still use a simplified/isotropic model for E's dynamics, our observer might be non-stationary (choosing among patrol trajectories) and our sensor might be only effective in an angular sector relative to O's current heading.

In Section \ref{s:time_dep_hjb}, we review a PDE approach to E's deterministic trajectory-planning, assuming the pointwise observability is fixed or predetermined.
Section \ref{s:se_games} describes the model and algorithms for SE games, taking adversarial planning into account.
Section \ref{s:numerics} focuses on numerical methods, 
while Section \ref{s:experiments} examines the Nash equilibrium policies for several test problems.
The limitations of our approach and directions for future work are discussed in Section \ref{s:conclusions}. 

\section{Evasive Path-Planning}
\label{s:time_dep_hjb}
We begin by considering a standard optimal control problem for the Evader under the assumption that the pointwise observability is fixed and fully known. 

E starts at some position $\bx \in \Omega$ at the time $t$, moves with a location-dependent speed $f(\bx) > 0$, and can instantaneously change its direction of motion.  
Its dynamics are given by:
\begin{equation}\label{eq:dynamics}
\by'(s) = f\bigl(\by(s)\bigr)\ba(s), \quad \by(t) = \bx,
\end{equation}
where $\ba : \mathbb{R} \to S^1$ is a measurable {\em control function} -- E's chosen (time-dependent) direction of motion.
Its path is constrained to stay within $\Omega \subset \mathbb{R}^2$ (\emph{e.g.,} avoiding any obstacles), and must reach the target $\xT$ by the time $T$.

Suppose that the pointwise observability function $K(\bx,s)$ is given (i.e., the observer's strategy is fixed).
Let $T_{\ba} = \min \{s\ge 0 ~ | ~ \by(s) = \xT \}$ be the time it takes the Evader to reach its target using control $\ba(\cdot)$.
E's goal is to minimize its cumulative observability:
\begin{equation}\label{eq:J} 
\J(\bx,t,\ba(\cdot)) = 
\begin{cases}
\int_t^{T_{\ba}} K(\by(s),s)ds, &T_{\ba} \le T \\
+\infty, &\text{ otherwise}
\end{cases}
\end{equation}
The value function $u(\bx,t)$ is then defined by
\begin{equation}\label{eq:valfunc}
u(\bx,t) = \inf_{\ba(\cdot)} \J(\bx,t,\ba(\cdot)).
\end{equation}

Standard arguments in optimal control theory (e.g.,  \cite[Chapters 3,4]{bardicapuzzodolcetta}) show that $u$ is the domain-constrained viscosity solution of a time-dependent Hamilton-Jacobi-Bellman (HJB) equation:
\begin{align}\begin{split}\label{eq:tdHJB}
\frac{\partial u}{\partial t} + \min_{|\ba|=1}\{f(\bx)\nabla u(\bx,t)\cdot \ba &+ K(\bx,t)\} = 0; \\
u\left( \bx, T\right) &= +\infty, \qquad \forall \x \neq \xT; \\
u\left(\xT,t\right) &= 0, \qquad \quad \forall t \leq T. 
\end{split}\end{align}
with an additional condition at $\partial\Omega \setminus \left\{ \x_T \right\}$ that the minimum is taken over the subset of control values that ensure staying inside $\overline{\Omega}$.
Wherever the value function $u$ is smooth, the optimal direction of motion is $\ba_* = -\nabla u / |\nabla u |$  and the above equation simplifies to a time-dependent Eikonal PDE
\begin{equation}\label{eq:tdEikonal}
\frac{\partial u}{\partial t} - f(\bx)\left|\nabla u(\bx,t) \right| + K(\bx,t) = 0.
\end{equation}

In general, \eqref{eq:tdEikonal} often does not have a classical solution, but always has a unique viscosity solution, coinciding with the value function of the optimal control problem \cite{bardicapuzzodolcetta}.
The optimal direction of motion $\ba_*$ is then uniquely determined almost everywhere on a set $\{(\x,t) \, | \, u(\x,t) < +\infty \}$ since $u$ is Lipschitz and thus differentiable except on a set of measure zero.
However, for those rare starting positions where $\nabla u$ is discontinuous, the optimal trajectory is not unique: there can be multiple optimal initial directions $\ba_*$, each of them producing a different optimal trajectory to $\xT.$  The numerical methods for solving this PDE and tracing the optimal trajectories will be covered in section \ref{s:numerics}.

\begin{figure}[htb]
\centering
\hspace*{-3mm}
$
\begin{array}[t]{cc}
\includegraphics[height=0.45\linewidth]{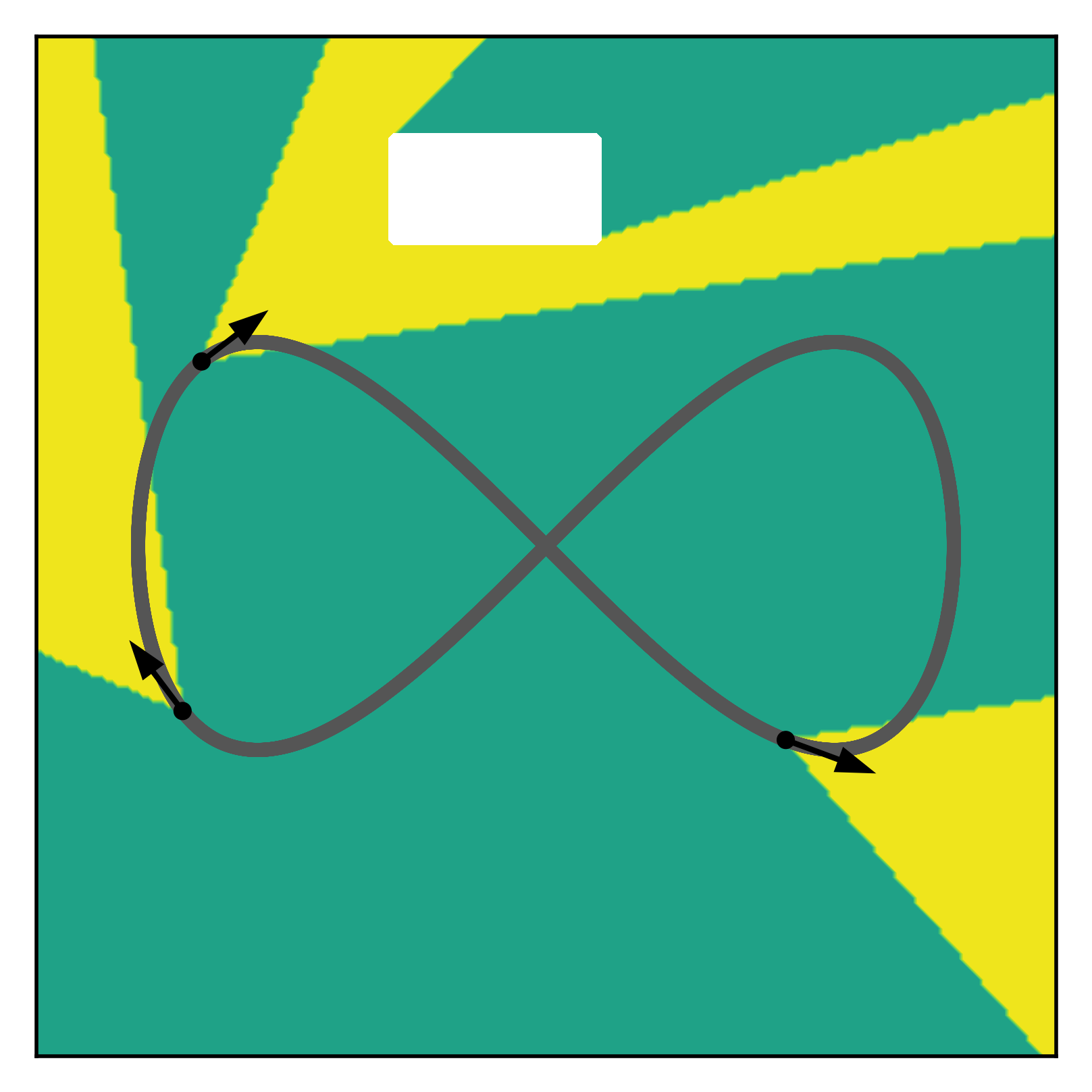} &
\includegraphics[height=0.45\linewidth]{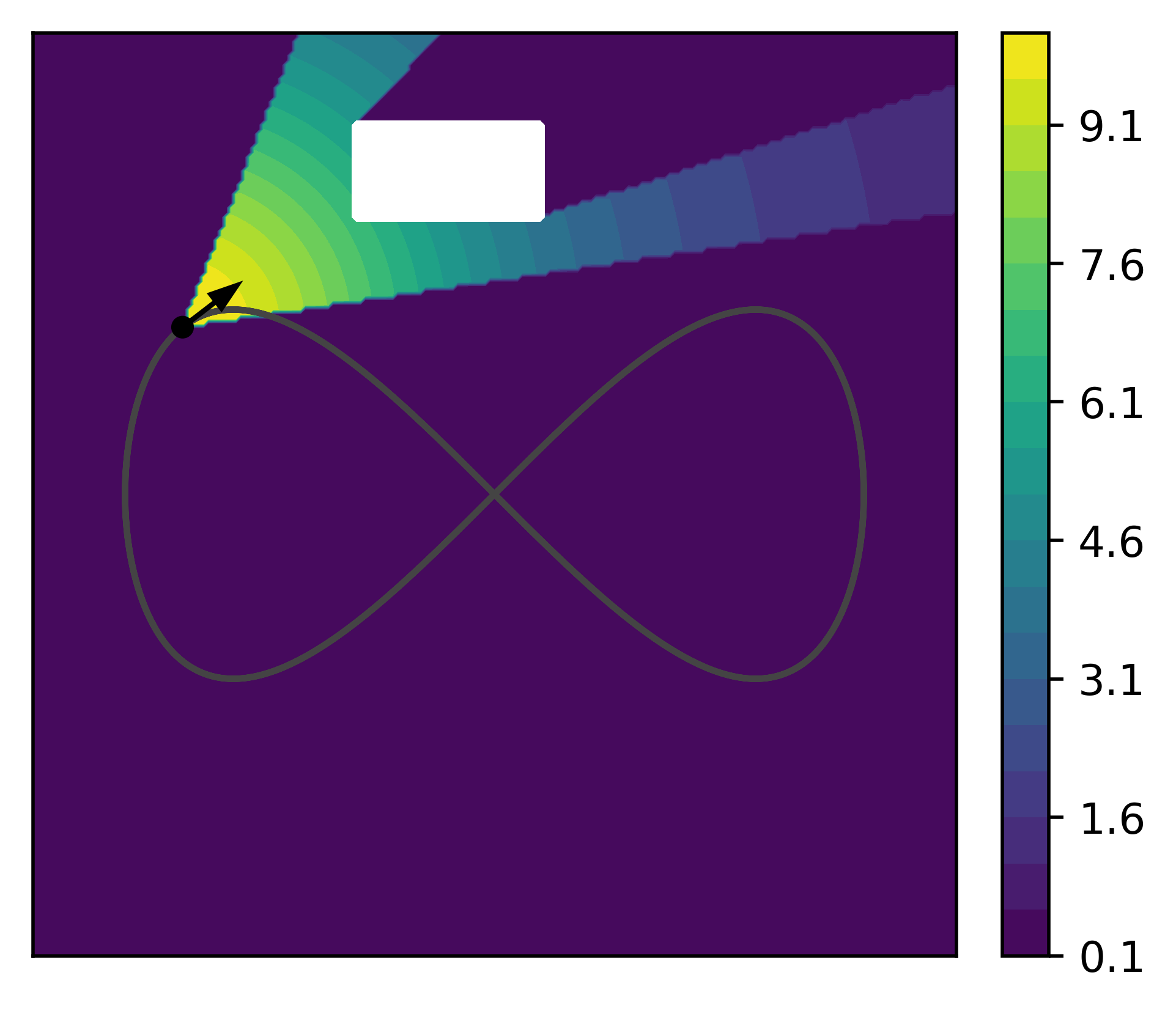}
\end{array}
$
\caption{\label{fig:obs}
LEFT: Visibility field (yellow) for 3 different Observer positions along a patrol trajectory (gray) in a domain with an occluding obstacle (white).   RIGHT: The contour plot of the pointwise observability $K(\x,t_1)$, corresponding to O's position $\z(t_1)$. 
}
\end{figure}

Before switching to a game-theoretic version of the problem, we remark on the definition of pointwise observability $K$.
Given O's current position $\z(t),$ one can define its {\em visibility field} $\V(t)$ restricted by obstacles and distance/angle limitations of the sensor.
Throughout this paper, we assume that all obstacles are fully occluding and the sensors are only effective in an angular sector of width $\alpha$ centered on observer's current heading $\bh(t).$
I.e., $\x \in \V(t)$ if the line segment $(\z(t), \x)$ stays in $\Omega$ without intersecting any obstacles and the observation angle (i.e., the angle that  $(\z(t), \x)$ makes with $\bh(t)$) is at most $\alpha/2$.
Within $\V(t)$, the pointwise observability should be a decreasing function of the distance $|\x-\z(t)|$. 
Assuming $\z(t)$ is known, we use
\begin{equation}\label{eq:obs}
K(\bx,t) = \
\begin{cases}
\frac{K_0}{\left|\bx-\bz(t)\right|^2 + 0.1} \, + \, \sigma, &\bx \in \V(t); \\
\sigma, &\bx \not\in \V(t);
\end{cases}
\end{equation}
where $K_0$ and $\sigma$ are positive constants.
This sector-restricted observability 
is illustrated  in Fig. \ref{fig:obs}.  

\section{Surveillance-Evasion Games}
\label{s:se_games}
Here we describe the zero-sum game between O and E, following the general framework developed in \cite{gilles2018surveillance}.
We assume that 
\begin{itemize}
    \item O can use any patrol trajectory from the set $\Z = \{\z_1(t), ... \z_r(t)\}$ known to both players;
    \item  E can use any of the infinitely many admissible trajectories that lead from $\xS$ to $\xT$ by the time $T$ while staying in $\Omega$ and avoiding the obstacles.
\end{itemize}
The payoff of our game is the (expected) cumulative observability, with E as a minimizer and O as a maximizer. 
Both players are assumed to make their decisions ahead of time and cannot change their mind based on any information gained while traveling through $\Omega.$

\subsection{The ``Evader-goes-second'' problem}
\begin{figure}[htb]
\centering
\hspace*{-3mm}
$
\begin{array}{cc}
\includegraphics[width=0.48\linewidth]{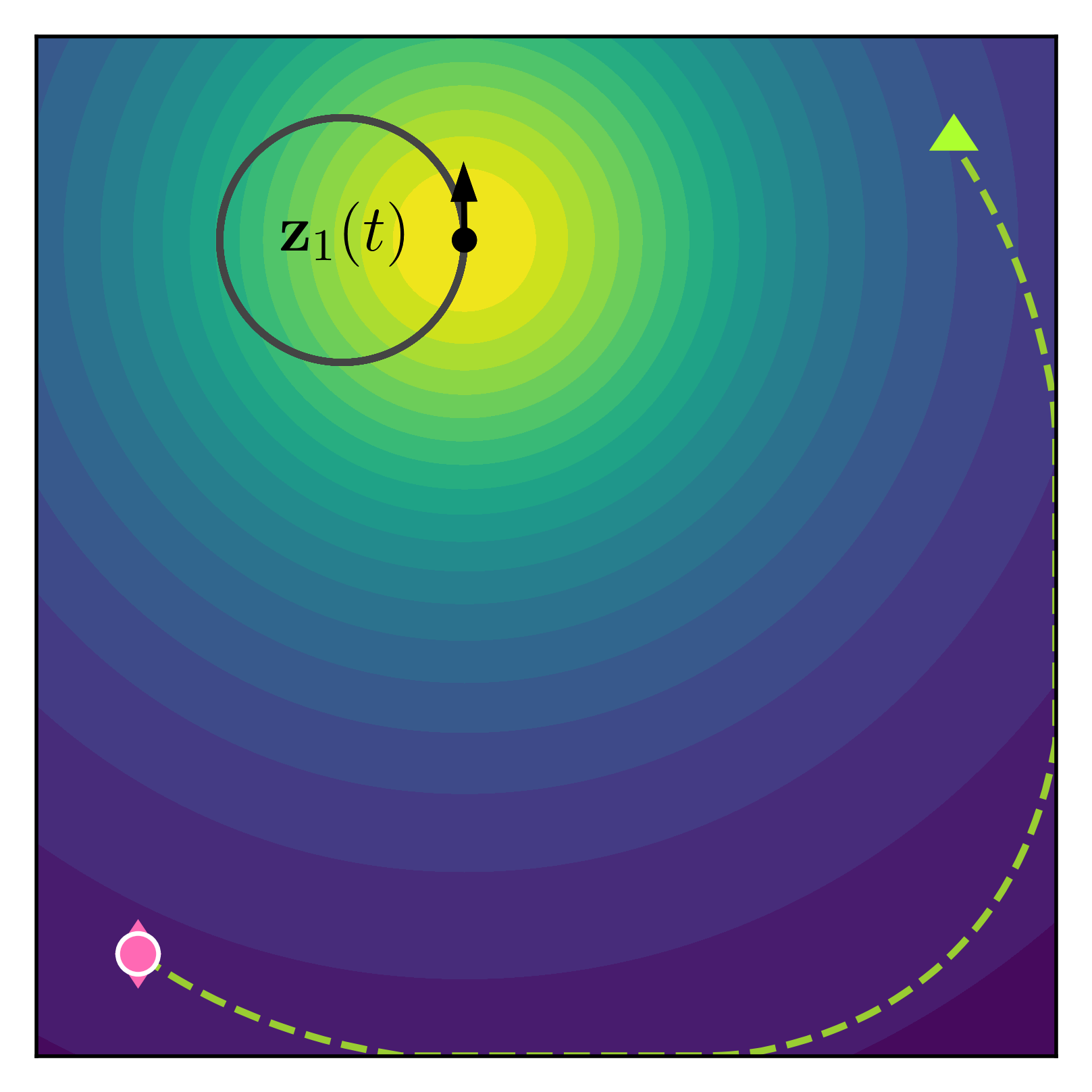} &
\includegraphics[width=0.48\linewidth]{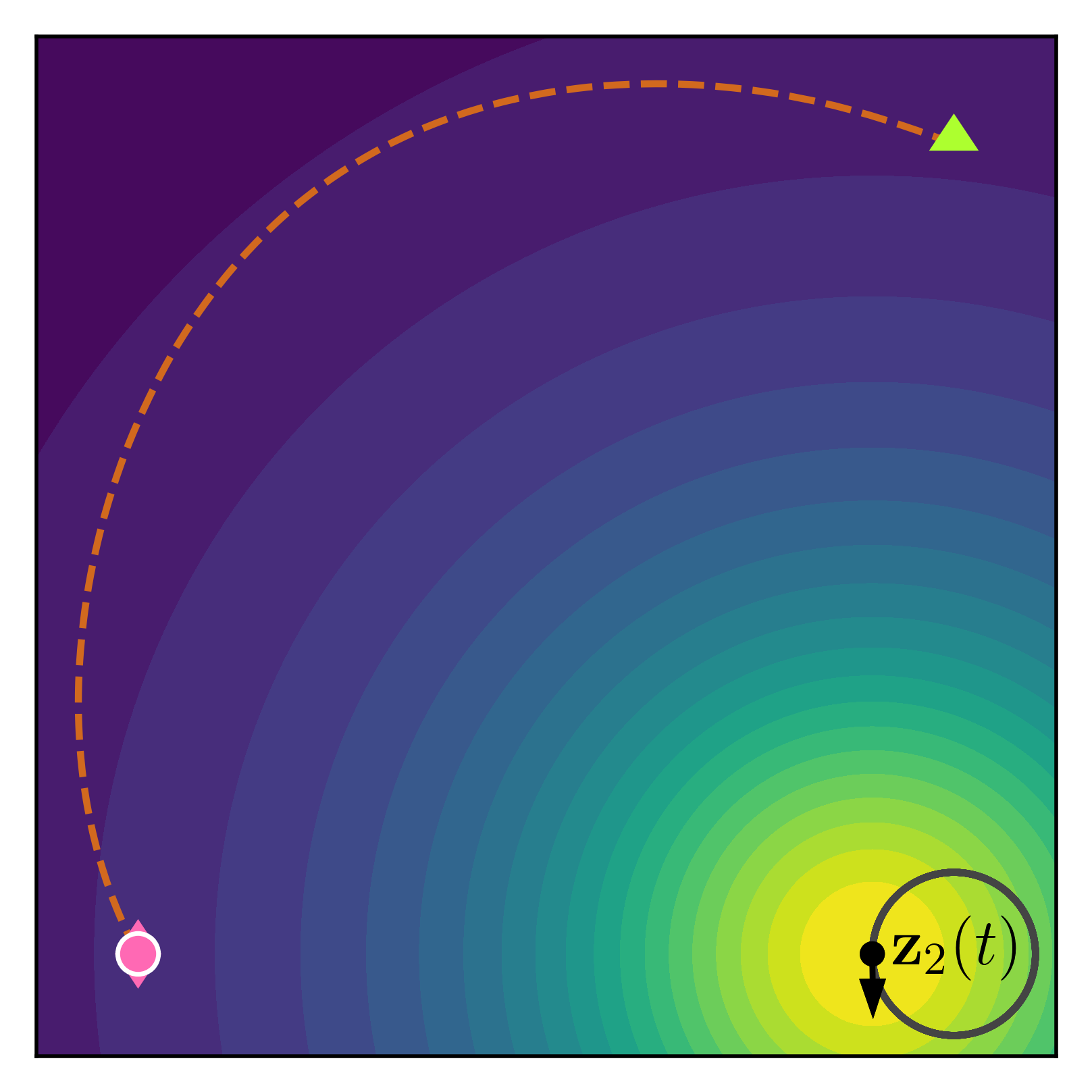}
\end{array}
$
\caption{\label{fig:pure}
E's optimal responses (shown in green and orange) corresponding to O's respective patrol trajectories.
The sensor is omnidirectional (i.e., $\alpha = 2\pi$).
The pointwise observability is time-dependent, and we show the contour plots for its {\em initial} version (i.e., $K_i(\x, 0)$) in both cases. 
}
\end{figure}

If O chooses the $i$-th patrol path, we can simply use $\z_i$ instead of $\z$ in \eqref{eq:obs} to
define the corresponding pointwise observability $K_i$, cumulative observability $\J_i$, and value function $u_i$.
The value of $u_i(\xS,0)$ yields the cumulative observability along the $\z_i$-optimal trajectory (i.e., E's best response to O's choice of $\z_i$.)
In Fig. \ref{fig:pure} we show two circular patrol trajectories (traversed by O with the same angular velocity) and E's optimal trajectory in response to each of them. If O were forced to make its decision deterministically before E, it would simply choose whichever $\z_i$ yields the highest
$u_i(\xS,0)$.

If O prefers to use a probability distribution $\lambda = (\lambda_1, \ldots, \lambda_r)$ over the set $\Z$ instead of choosing a specific $\z_i$, we can similarly define the expected pointwise observability $K^{\lambda} = \sum_{i=1}^r \lambda_i K_i$ and the corresponding expected cumulative observability $\J^{\lambda}=\sum\limits_{i=1}^r \lambda_i \J_i$.  To find E's best response, the corresponding value function $u^{\lambda}$ can be computed by solving \eqref{eq:tdEikonal} with $K^{\lambda}$ used in place of $K.$   We note that $u^{\lambda}$ is {\em not} a linear combination of $u_i$'s since this PDE is nonlinear.

\begin{figure}[htb]
\centering
\hspace*{-3mm}
$
\begin{array}{cc}
\includegraphics[width=0.48\linewidth]{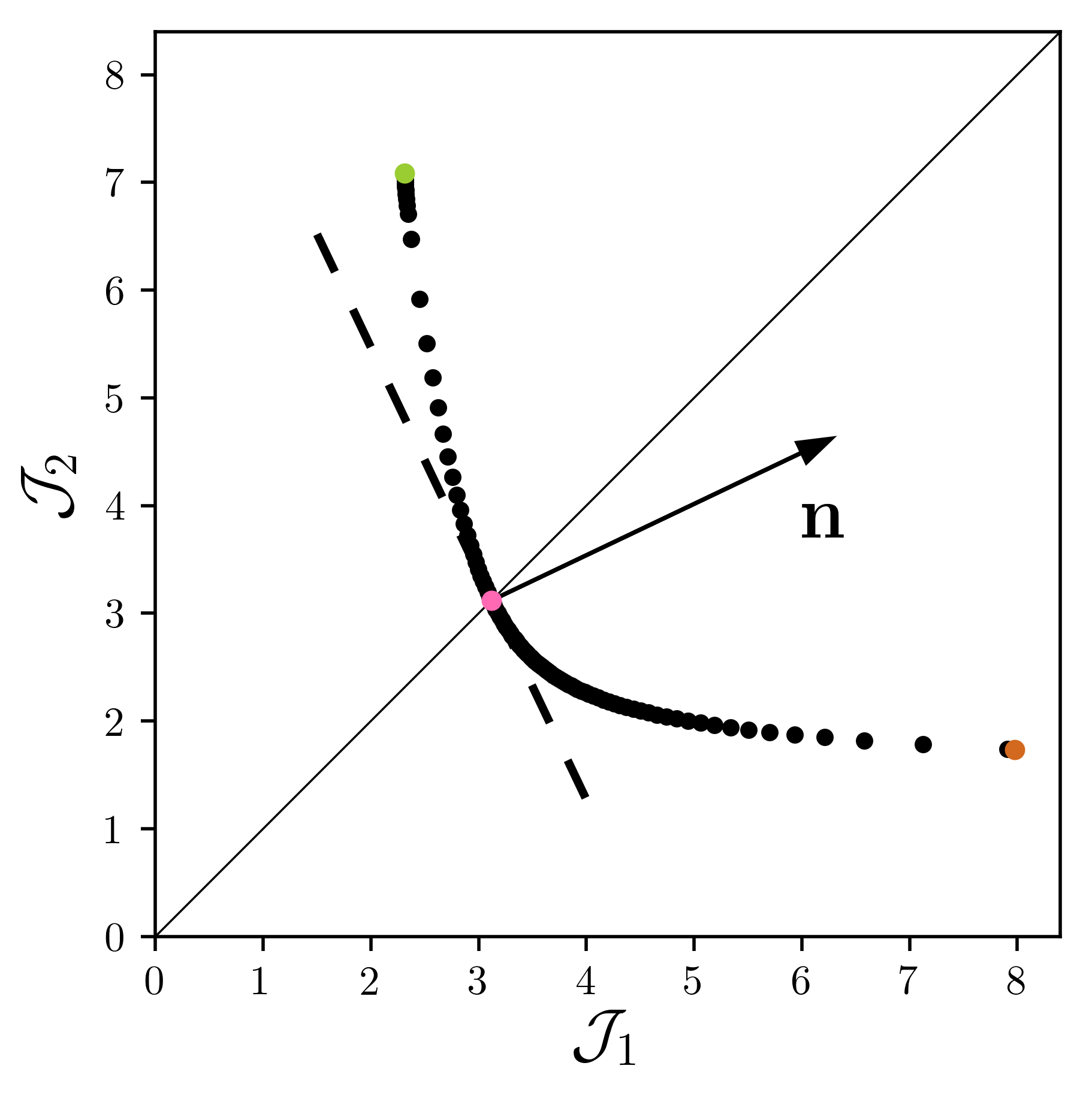} &
\includegraphics[width=0.48\linewidth]{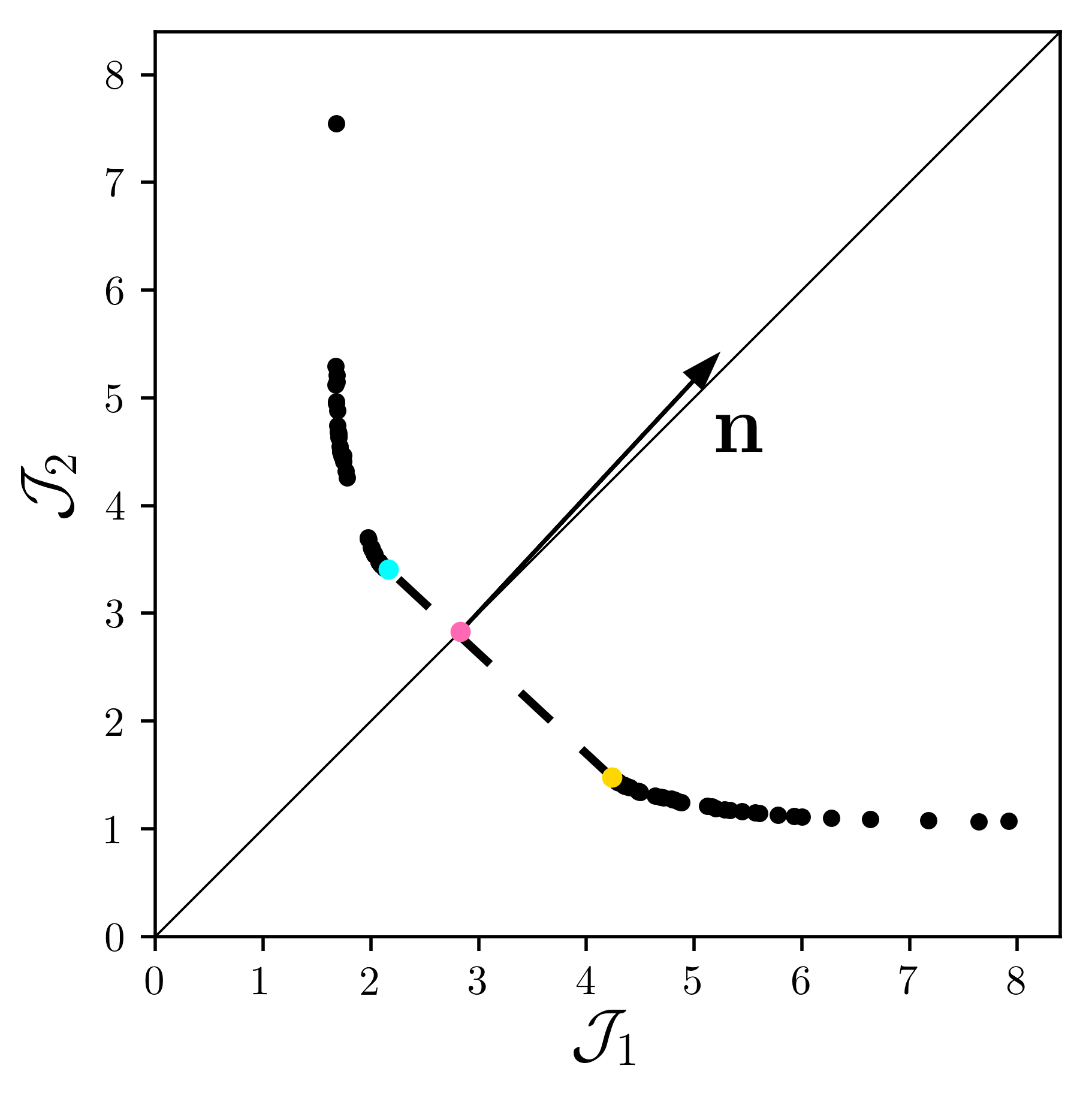}\\
(A) & (B)
\end{array}
$
\caption{\label{fig:PF}
The convex part of Pareto Fronts (PFs) corresponding to two examples from section \ref{s:experiments}.
PF approximations generated by using 101 $\lambda$ values, solving for $u^{\lambda}$, and computing $(\J_1, \J_2)$ for $\lambda$-optimal trajectories.  The dashed line shows a support hyperplane to PF, with $\boldsymbol{n}$ showing a scaled version of $\lambda_*$.    
Evader's 
optimal policy (shown in magenta) corresponds to the intersection of that hyperplane with a ``central ray'' $\J_1=\J_2$. \newline
LEFT: Example 1 has only one $\lambda_*$-optimal trajectory and E's optimal policy is deterministic; see Fig. \ref{fig:ex1}. $\J_1$-optimal and $\J_2$-optimal trajectories are also represented on PF by the green and orange markers respectively. \newline
RIGHT: In Example 2, E's optimal policy assigns non-zero probabilities to both $\lambda_*$-optimal trajectories (cyan and yellow); see Fig. \ref{fig:ex2}.
}
\end{figure}

This approach is also related to methods for {\em multi-objective} optimal control. It is generally impossible to minimize the observability with respect to all patrol trajectories simultaneously.  Instead, one can consider the task of finding {\em Pareto-optimal} controls.  A control $\tilde{\ba}(\cdot)$ dominates $\ba(\cdot)$ if 
$\J_i \left(\xS,0,\tilde{\ba}(\cdot)\right) \, \leq \, 
\J_i \left(\xS,0,\ba(\cdot)\right)$ for all $i$, with the inequality strict for at least one $i$.  We will say that $\ba(\cdot)$ is Pareto-optimal if it is not dominated by any other control.  Plotting the point $(\J_1, \ldots, \J_r)$ for each Pareto-optimal control we obtain a Pareto Front (PF), illustrated in Fig. \ref{fig:PF} for two examples from section \ref{s:experiments}.
Each vector $\lambda$ is normal to a support hyperplane of PF (on which $\J^{\lambda} = u^{\lambda}(\xS,0)$).  If all the elements of $\lambda$ are positive, it is easy to show that any $\lambda$-optimal control is also Pareto-optimal. 
This is the basis of a {\em weighted-sum scalarization} approach to general multi-objective trajectory planning introduced in \cite{MitchellSastry}. 
Unfortunately, it approximates only convex parts of PF \cite{das1997closer}.  Moreover, it requires imposing a fine grid on the space of $\lambda$'s and solving the PDE for each $\lambda$-gridpoint. 
Alternative methods can recover the entire PF \cite{KumarVlad,UTRC_2}, but are more computationally expensive.  
Luckily, for SE games we only need the points on PF corresponding to the players' {\em mutually} optimal policies; below we explain how this can be done through scalarization with a much smaller number of PDE solves.

\noindent
O's goal is to find a $\lambda$ maximizing
$$
G(\lambda) 
\; = \; \min\limits_{\ba(\cdot)} 
\sum\limits_{i=1}^r 
\lambda_i \J_i \left(\xS,0,\ba(\cdot)\right)
\; = \; 
u^{\lambda}(\xS,0).
$$
Since $G$ is a pointwise minimum of functions linear in $\lambda$, we know that $G(\lambda)$ is concave.  Thus, O's optimal pdf $\lambda_*$ can be found by standard methods of convex optimization. Our implementation relies on the projected supergradient ascent method \cite[Chap.~8]{beck2017first} to maximize $G$ on the $r$-dimensional simplex of probability densities. During the iterative improvement of $\lambda$, we use gradient descent in $u^{\lambda}$ starting from $(\xS,0)$ to find E's $\lambda$-optimal trajectory.  Integrating each of $K_i$'s along that trajectory we obtain the supergradient of $G(\lambda)$, which is then used to update the $\lambda$ for the next iteration.
See Algorithm 4.2 in \cite{gilles2018surveillance}.

\subsection{Nash Equilibria}
It is restrictive to assume that E is always reacting to O's decisions.  To get a more general interpretation of the game's value, we need to search for a {\em Nash equilibrium}.
E can also define its policy probabilistically, by choosing a probability distribution $\theta$  
over the infinite set of admissible controls.  
Given the players' policies, the expected payoff is
$$
P(\lambda, \theta) 
\! = \!
\E_{\lambda, \theta}\biggl[ \J\!\left(\xS,0,\ba(\cdot)\right)\biggr] \! =\! \sum\limits_{i=1}^r \lambda_i \E_{\theta} \biggl[\J_i\!\left(\xS,0,\ba(\cdot)\right)\biggr],
$$
where $\ba(\cdot)$ is a random admissible control chosen according to $\theta$.
A pair of policies $(\lambda_*, \theta_*)$ forms a Nash equilibrium if neither player can improve the payoff by changing its policy unilaterally.  I.e.,
$$
P(\lambda, \theta_*) \leq 
P(\lambda_*, \theta_*) \leq 
P(\lambda_*, \theta), \quad \forall \lambda, \theta. 
$$
There always exists at least one Nash equilibrium in the set of mixed/probabilistic policies.  Moreover, in zero-sum two-player games, the payoff is the same for all Nash equilibria and is thus used to define the {\em game's value}; e.g., see \cite{osborne1994course}.  Thus, when each player uses a Nash policy, it no longer matters which of them ``goes first.''

By the min-max theorem, the value of this game can be found by computing
\begin{equation}\label{eq:maxmin}
P(\lambda_*, \theta_*) \; = \;
\min\limits_{\theta} \max\limits_{\lambda}  
P(\lambda, \theta) \; = \; \max\limits_{\lambda} \min\limits_{\theta} P(\lambda, \theta).
\end{equation}
In finite games, this is typically accomplished by linear programming \cite{osborne1994course}, but since E has infinitely many possible paths, our game is {\em semi-infinite} \cite{tijs1976semi,raghavan1994zero} and this approach is inapplicable.  Luckily, it is easy to see that the last expression in \eqref{eq:maxmin} is equal to $\max_\lambda G(\lambda)$ and the previous subsection explains how to find the maximizer $\lambda_*$ efficiently.

Once the algorithm converges to O's optimal $\lambda_*$, we need to find the second half of the Nash equilibrium: E's optimal probability distribution $\theta_*$ over the set of admissible trajectories.  We accomplish this using the following properties, proven for a stationary Observer in \cite{gilles2018surveillance}, which similarly hold for the time-dependent case considered here.
Suppose $(\lambda_*, \theta_*)$ is a Nash equilibrium, $\I = \{i | \lambda_{*,i} > 0 \}$, and $\A$ is the set of all
controls that have positive probability under $\theta_*.$
\begin{enumerate}
\item
If $\lambda_*$ maximizes $G(\lambda)$, then there always exists a probability distribution $\theta_*$ over Evader's Pareto-optimal trajectories such that
$(\lambda_*, \theta_*)$ is a Nash equilibrium and the set $\A$ has at most $|\I| \leq r$ elements.
\item 
If $\ba(\cdot) \in \A$ then $\ba(\cdot)$ is $\lambda_*$-optimal.\\  I.e., 
$\sum\limits_{i=1}^r 
\lambda_i  
\J_i \left(\xS,0,\ba(\cdot)\right)
\, = \, u^{\lambda_*} (\xS,0).$
\item
If $i \in \I$, then $\E_{\theta_*} \left[ 
\J_i \left(\xS,0,\ba(\cdot)\right)
\right] \, = \, u^{\lambda_*} (\xS,0).$
\end{enumerate}
For all starting positions where $u^{\lambda_*}$ is differentiable, the optimal direction of motion is $ -\nabla u^{\lambda_*} / |\nabla u^{\lambda_*}|$ and the optimal control is unique. Thus, it is natural to expect that, {\em generically}, there should be only one $\lambda_*$-optimal control $\ba_*(\cdot)$ leading from $\xS$ to $\xT$ by the deadline $T$. This would imply that $\A = \left\{ \ba_*(\cdot) \right\}$; i.e., E's half of the Nash Equilibrium is pure/deterministic and
$\J_i \left(\xS,0,\ba_*(\cdot)\right) = G(\lambda_*) = u^{\lambda_*} (\xS,0)$
for all $i \in \I$.  (We say that $\ba_*(\cdot)$ is the intersection of the Pareto Front and the ``central ray'' in the $|\I|$-dimensional cost space.)
This situation, illustrated in Fig. \ref{fig:PF}(A) is indeed quite common.  But surprisingly, it is not generic; i.e., for many problems $\A$ must include more than one $\lambda_*$-optimal control.  This is a counter-intuitive side-effect of changing $\lambda$ to maximize $G$: $K^{\lambda}$ keeps changing until $(\xS,0)$ falls onto a ``shockline'', a set of discontinuities of $\nabla u^{\lambda_*}.$  
A numerically implemented gradient descent in $u^{\lambda_*}$ usually yields only one $\lambda_*$-optimal trajectory.  Approximating {\em all} of them is much harder.  We accomplish this by perturbing $\lambda_*$ in various directions and computing the resulting unique optimal trajectories corresponding to the perturbed $(\lambda_* + \delta \lambda)$.  We refer readers to Alg. 4.3 in \cite{gilles2018surveillance} for a detailed discussion.

\section{Numerical Implementation}
\label{s:numerics}

\subsection{Numerics for HJB equation}
\label{ss:hjb_numerics}
\begin{figure}[h]
\begin{center}
\begin{tikzpicture}
    \draw (0,0) -- (0,1) -- (1,1) -- cycle;
    \draw (0,0) -- (0,1) -- (-1,1) -- cycle;
    \draw (0,0) -- (1,0) -- (1,1) -- cycle;
    \draw (0,0) -- (1,0) -- (1,-1) -- cycle;
    \draw (0,0) -- (0,-1) -- (-1,-1) -- cycle;
    \draw (0,0) -- (0,-1) -- (1,-1) -- cycle;
    \draw (0,0) -- (-1,0) -- (-1,1) -- cycle;
    \draw (0,0) -- (-1,0) -- (-1,-1) -- cycle;
    \draw (0,1.2) node {(i,j+1)};
    \draw (1.35,1.2) node {(i+1,j+1)};
    \draw (1.5,0) node {(i+1,j)};
    \draw (1.35,-1.25) node {(i+1,j-1)};
    \draw (0,-1.25) node {(i,j-1)};
    \draw (-1.35,-1.25) node {(i-1,j-1)};
    \draw (-1.5,0) node {(i-1,j)};
    \draw (-1.35,1.2) node {(i-1,j+1)};
    \draw (0.25,0.6) node {2};
    \draw (0.6,0.25) node {1};
    \draw (0.6,-0.25) node {8};
    \draw (0.25,-0.6) node {7};
    \draw (-0.25,-0.6) node {6};
    \draw (-0.6,-0.25) node {5};
    \draw (-0.6,0.25) node {4};
    \draw (-0.25, 0.6) node {3};
\end{tikzpicture}
\end{center}
\caption{Nine-point stencil used to discretize HJB equation at $(\x_{i,j}, t_k).$}
\label{fig:8pt_stencil}
\end{figure}
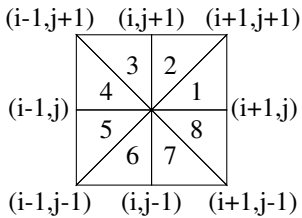

We numerically solve \eqref{eq:tdEikonal} using a time-explicit first-order upwind finite-difference scheme on a 9-point stencil.
Let $\x_{i,j}=(i h, jh)$, $t_k=k\Delta t$, and 
$U^k_{i,j} \approx u\left(\x_{i,j},t_k\right)$.
We will solve the equation backwards in time, computing $U^{k-1}_{i,j}$ using the values of $U$ at time slice $k$. 
More specifically, we will compute an update from each of the eight triangular simplices generated by $\x_{i,j}$ and its eight nearest neighbors (see Fig. \ref{fig:8pt_stencil}). 

We first compute the first-order approximation $(U_x,U_y) \approx \nabla u(\x_{i,j}, t_k)$ based on a specific simplex. 
We then check the {\em upwinding condition}, requiring the approximate gradient to point from that same simplex.
E.g., for simplex 1
\begin{equation*}
    U_x = \frac{U^k_{i+1,j} - U^k_{i,j}}{h}, \qquad U_y = \frac{U^k_{i+1,j+1}-U^k_{i+1,j}}{h},
\end{equation*}
and the corresponding upwinding condition is $U_x \le U_y \le 0$.

Focusing on a single simplex $l\in\{1,...,8\}$ with vertices $\x_{i,j}, \x^{l,1},$ and $\x^{l,2}$, the update from that simplex is
\begin{equation*}
    U_l = U^k_{i,j} - \Delta t f(\x_{i,j}) D^k_{i,j} + \Delta t K(\x_{i,j}, t_k),
\end{equation*}
where $D^k_{i,j}$ approximates $|\nabla u (\x_{i,j}, t_k)|.$
If the upwinding condition is satisfied, we use 
$D^k_{i,j} = \sqrt{U^2_x + U^2_y}.$  Otherwise, we resort to a ``one-sided'' semi-Lagrangian update, with
$$
D^k_{i,j} \; = \; \max \left(
\frac{U^k_{i,j} - U(\x^{l,1}, t_k)}{|\x_{i,j} - \x^{l,1}|}, \,
\frac{U^k_{i,j} - U(\x^{l,2}, t_k)}{|\x_{i,j} - \x^{l,2}|}, \,
0
\right).
$$

We then set $U^{k-1}_{i,j}$ to be the smallest of simplex-specific updates $\left\{U_l\right\}$.
The boundary conditions and terminal conditions are enforced by setting
\begin{equation*}
U^k_{i,j} = 
\begin{cases}
0, & \quad \x_{i,j} = \xT \text{ and } t_k \leq T;\\
+\infty, & \quad \x_{i,j} \not\in\overline{\Omega} 
\, \text{ \bf or } \, (\x_{i,j} \not=\xT \text{ and } t_k=T).
\end{cases}
\end{equation*}

As long as the CFL stability condition $\Delta t \le h / (\max f(\x_{i,j}))$ is satisfied, this method is monotone and consistent \cite{Falcone2013}, and therefore converges to the viscosity solution of \eqref{eq:tdEikonal} under grid refinement \cite{barles1991convergence}.

\subsection{Time-dependent Visibility}
\label{ss:time_dep_shadows}
Our definition of $\V_i(t)$ requires computing the set of points that have direct line of sight from $\bz_i(t)$. 
Computational efficiency is very important as this needs to be computed for every time slice.
Following \cite{tsai2004visibility}, we first construct a continuous function $\phi(\bx)$ such that
\begin{equation*}
\phi(\bx)
\begin{cases}
> 0, &\bx \in \Omega, \\
= 0, &\bx \in \partial\Omega, \\
< 0, &\bx \in \mathbb{R}^2 \setminus \overline{\Omega} 
\quad \text{(e.g., inside obstacles).} 
\end{cases}
\end{equation*}

Within each time slice $t$, we numerically solve the quasi-variational inequality
\begin{align*}
\max\{\nabla \psi^t(\bx) \cdot \br(\bx), \; \psi^t(\bx) - \phi(\bx)\}  &= 0,\\
\psi^t(\bz(t)) &= \phi(\bz(t)),
\end{align*}
where $\br(\bx)$ is a unit vector pointing from $\bz(t)$ to $\bx.$
Then $\{\psi^t \ge 0\}$ will define the set of points that are not occluded by obstacles.
After computing $\psi^t$, we then enforce any additional distance or angular restrictions to get $\V_i(t)$.
Combined with \eqref{eq:obs}, this allows us to pre-compute and store $K_i(\bx,t)$ for each patrol trajectory $\z_i(t)$ before starting to plan E's optimal responses to various Observer policies $\lambda$.

\subsection{Time-dependent Optimal Path Tracer}
\label{ss:path_tracer}
Given our numerical solution $U^k_{i,j}$, we now wish to recover the Evader's optimal trajectory $\by(t)$.
We approximate the trajectory with a series of points $\by^m \approx \by(m\Delta t)$, where $\Delta t$ is the timestep used to numerically solve \eqref{eq:tdEikonal}.

E's initial position is $\by^0 = \by(0) = \xS$.
The rest of the path can be found in a semi-Lagrangian manner by computing
\begin{align*}
\ba^m_* &= \argmin_{|\ba| \le 1} \bigl\{ \widetilde{U}\left(\by^m + \Delta t f(\by^m) \ba, \, m\Delta t \right) \bigr\}, \\
\by^{m+1} &= \by^m + \Delta t \ba^m_*. 
\end{align*}
where $\widetilde{U}$ is the function 
obtained from the grid values of $U$ by trilinear interpolation.
We terminate the path once we get close enough to the target so that $\left| \by_m - \xT \right| \le f(\by_m)\Delta t$.
In our numerical experiments, we optimize 
over $\{|\ba|=1\} \cup \{(0,0)\}$ instead of the entire unit disk $\{|\ba|\le 1\}$.

\section{Numerical Experiments}
\label{s:experiments}
\begin{figure*}[p!]
\centering
$
\arraycolsep=1pt\def\arraystretch{0.1}
\begin{array}{cccc}
\includegraphics[height=0.18\textheight]{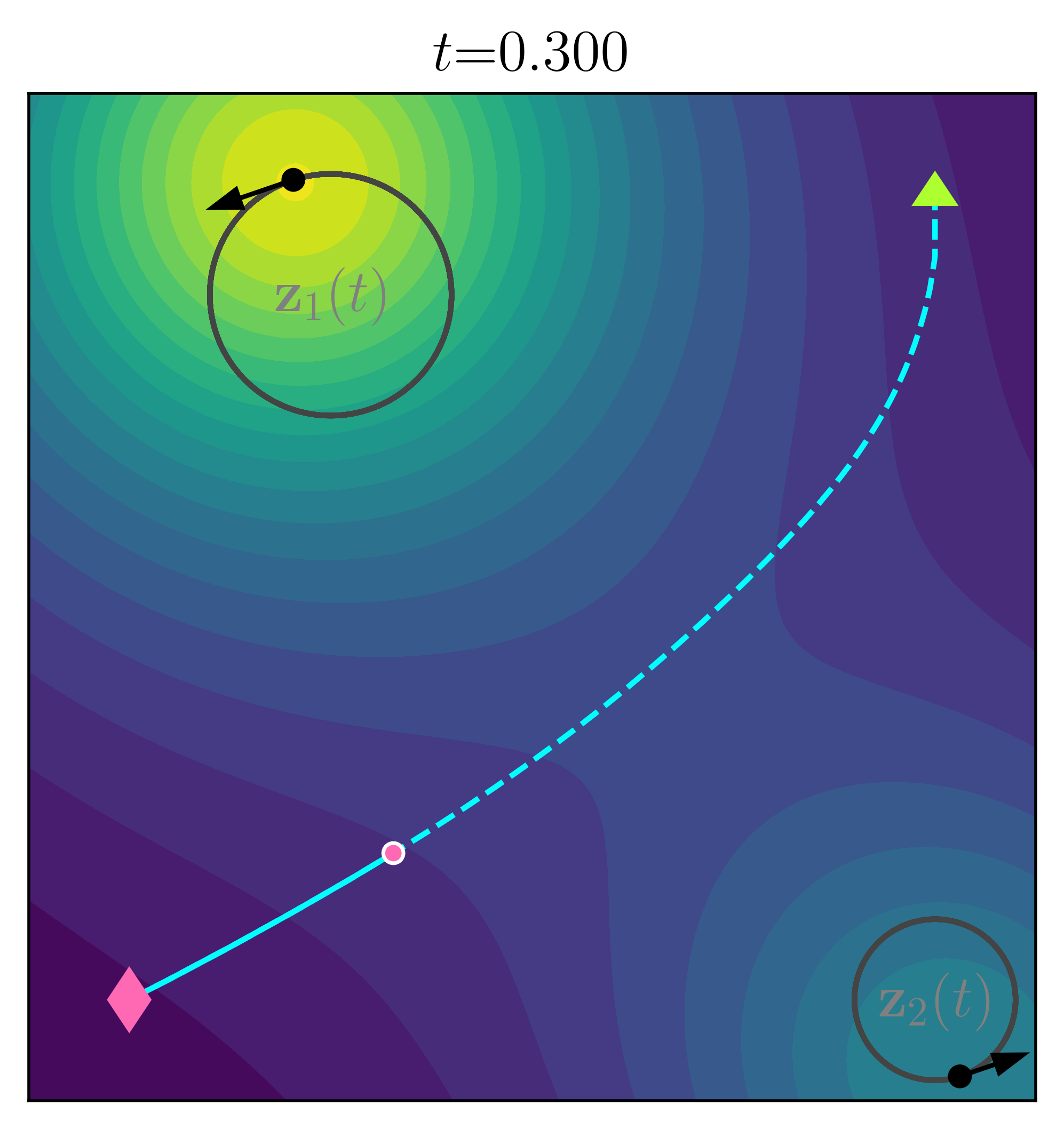} &
\includegraphics[height=0.18\textheight]{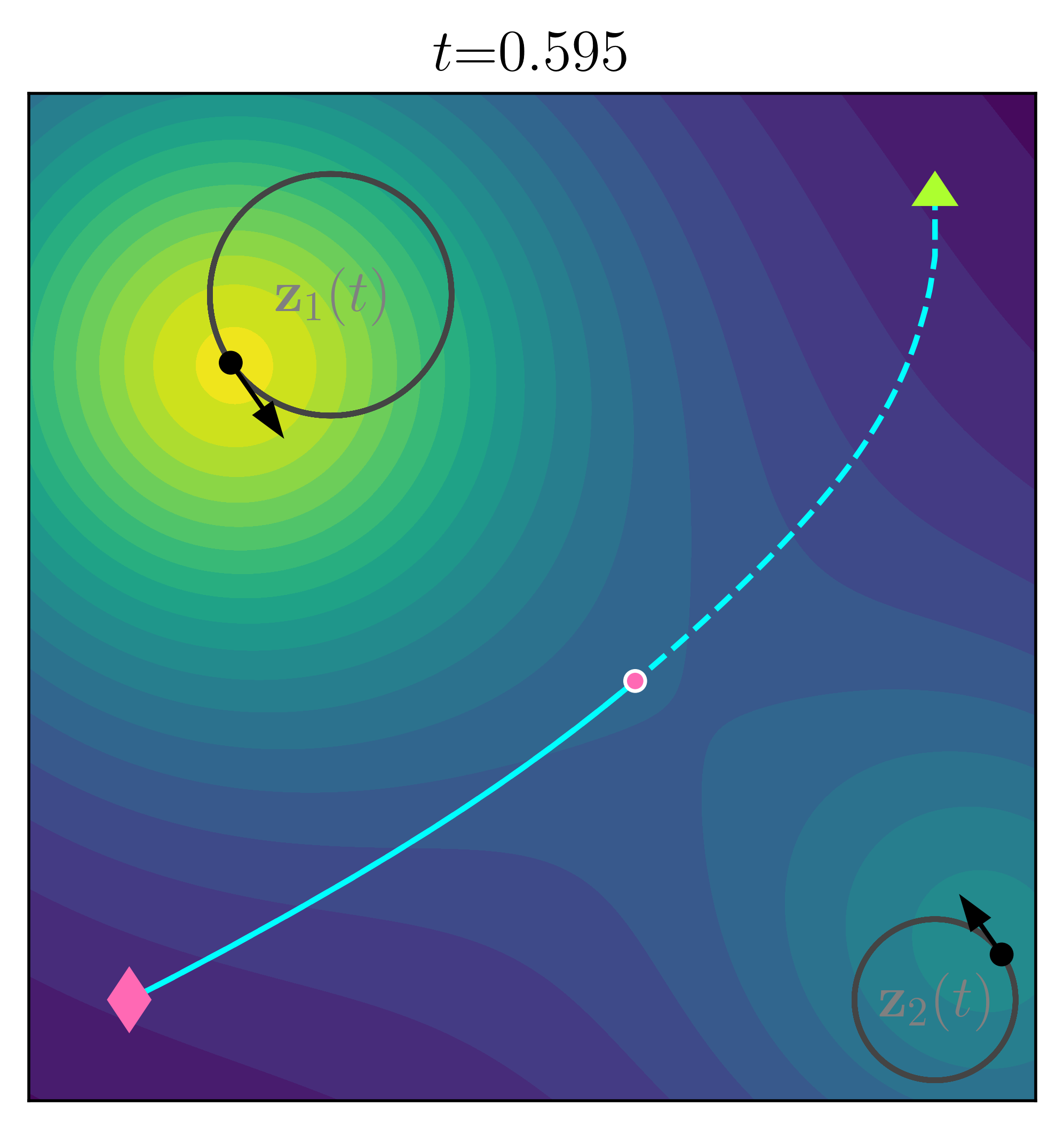} &
\includegraphics[height=0.18\textheight]{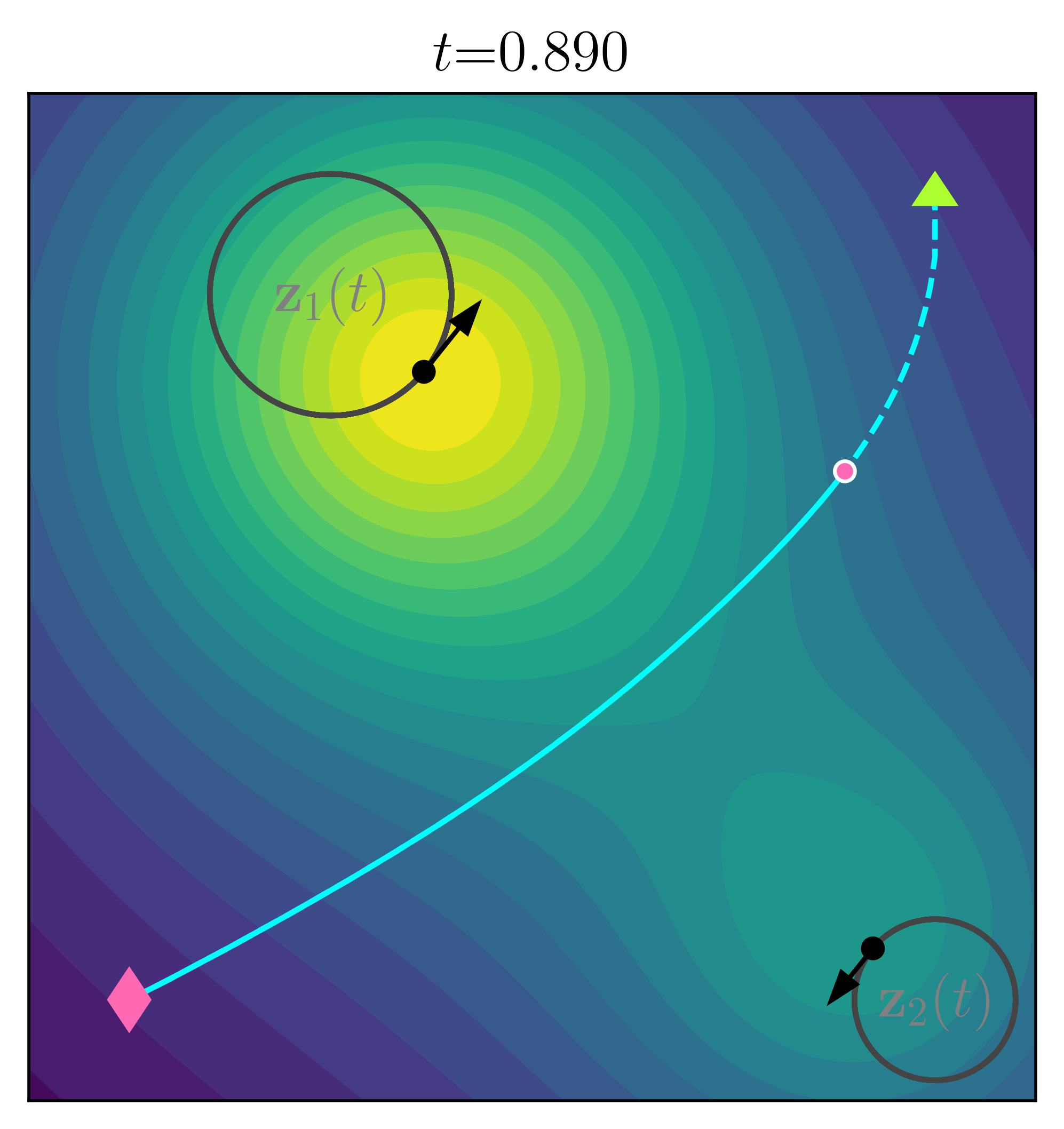} &
\includegraphics[height=0.18\textheight]{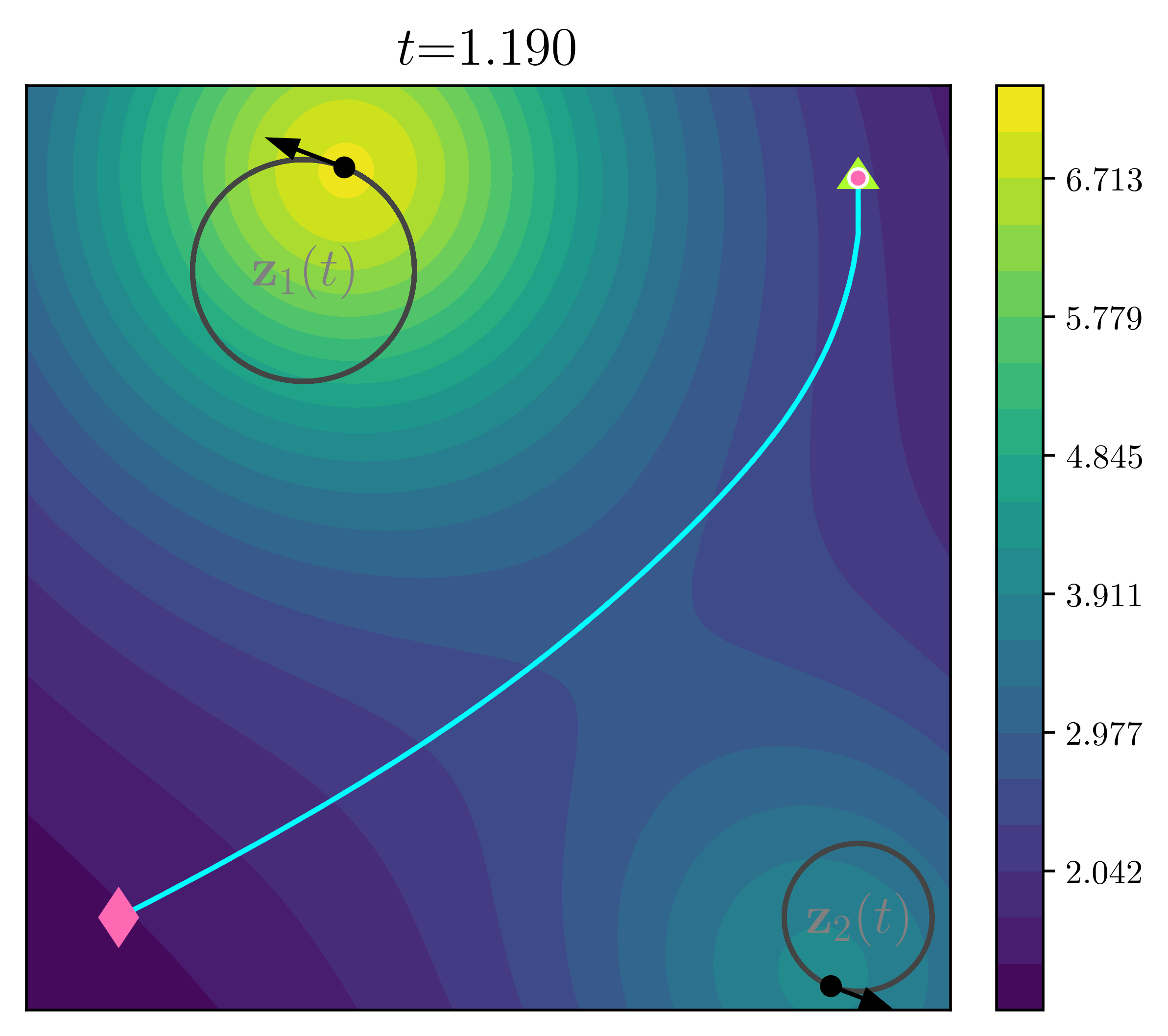} 
\end{array}
$
\caption{Four snapshots of the Nash Equilibrium solution for Example 1.  Observer's optimal policy is $\lambda_*=(0.67,0.33).$  Observer's patrol trajectories (in gray) and Evader's unique $\lambda_*$-optimal trajectory (magenta) shown on top of the contour plots of $K^{\lambda_*}(\x,t).$  Observer's possible position on each patrol trajectory is shown by a black circle, with an arrow indicating Observer's current heading.  Evader's starting, current, and terminal positions are shown by a magenta diamond, a magenta circle, and a green triangle respectively. \label{fig:ex1}}
\end{figure*}

\begin{figure*}[p!]
\centering
$
\arraycolsep=1pt\def\arraystretch{0.1}
\begin{array}{cccc}
\includegraphics[height=0.18\textheight]{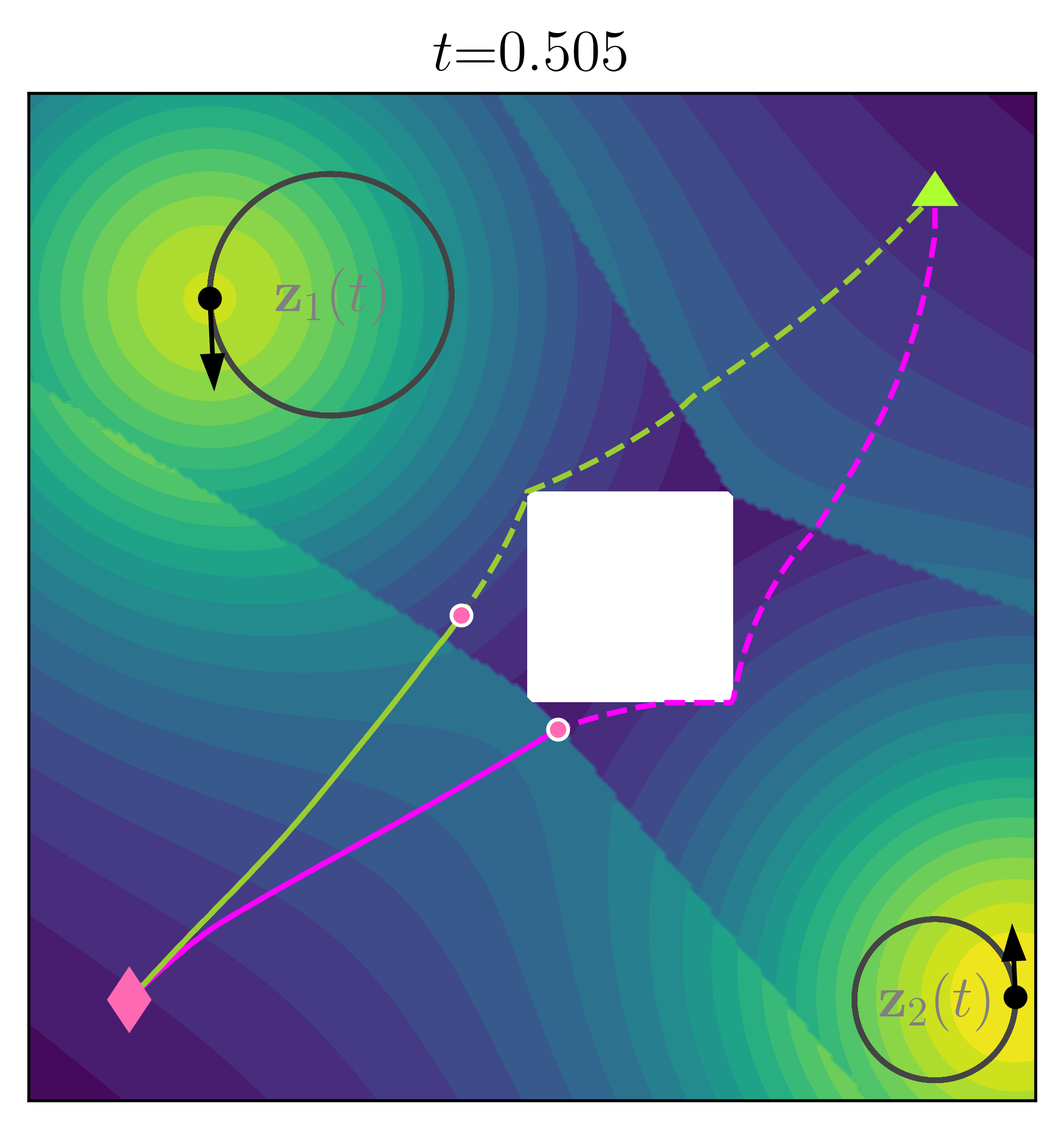} &
\includegraphics[height=0.18\textheight]{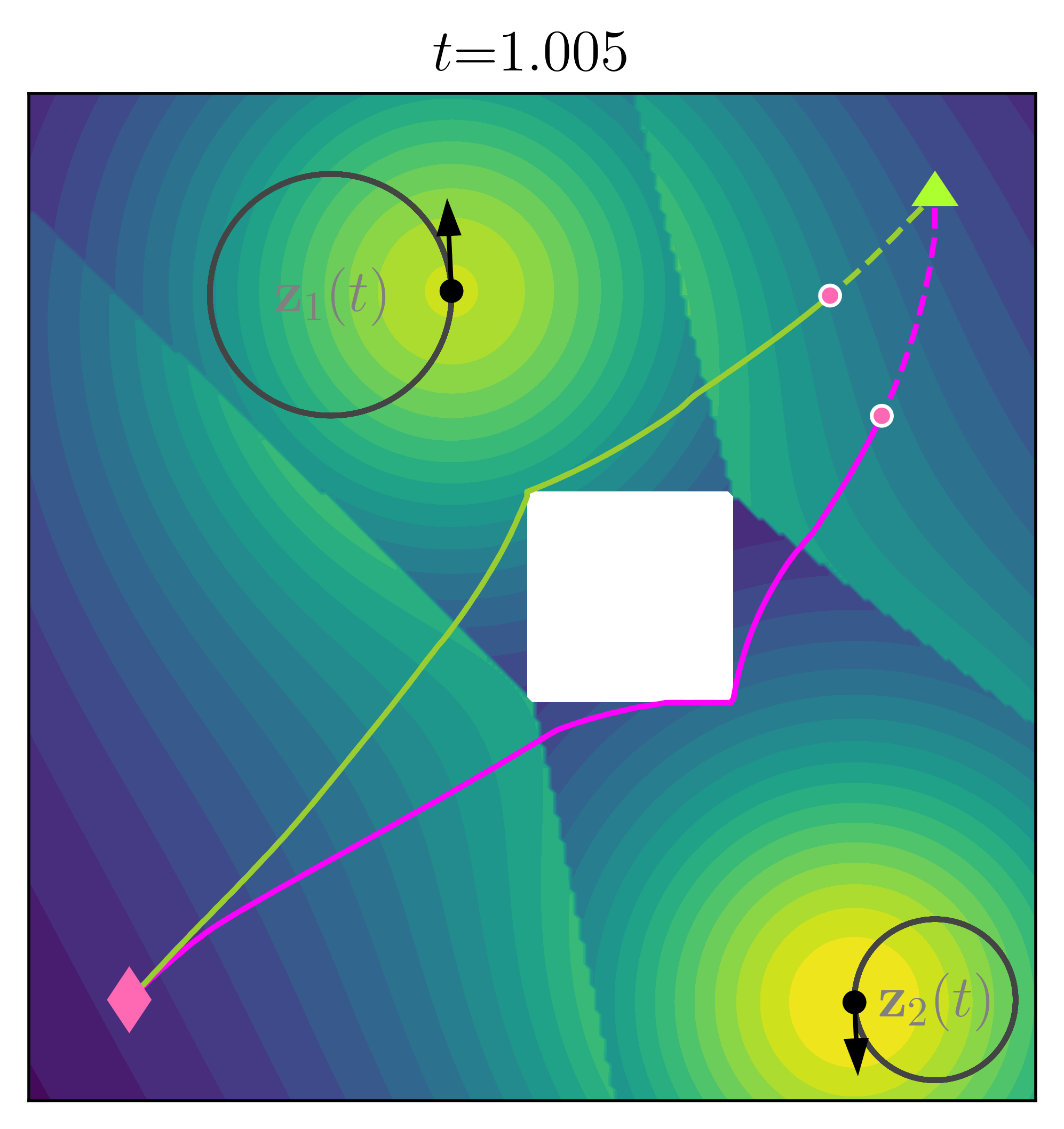} &
\includegraphics[height=0.18\textheight]{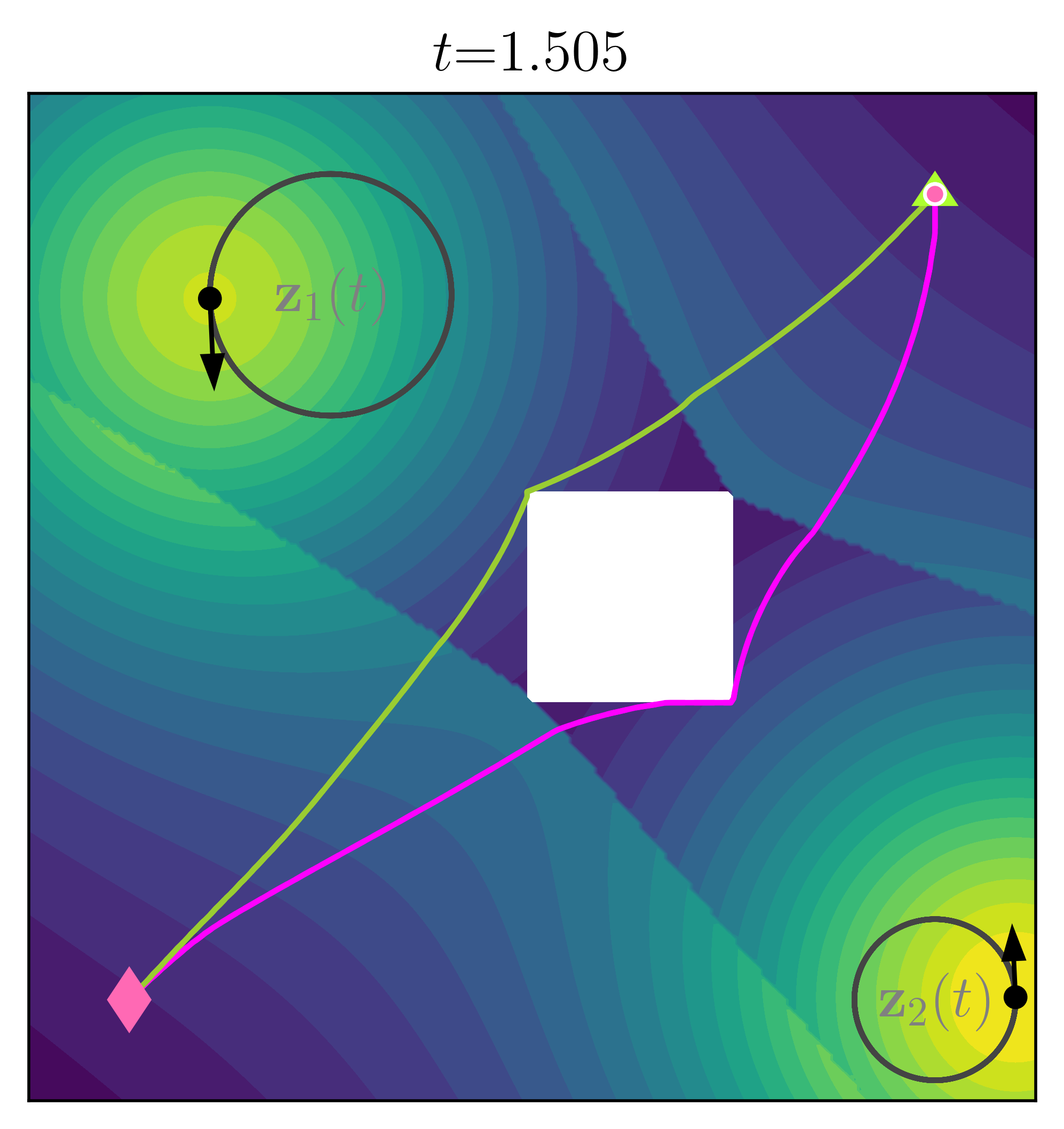} &
\includegraphics[height=0.18\textheight]{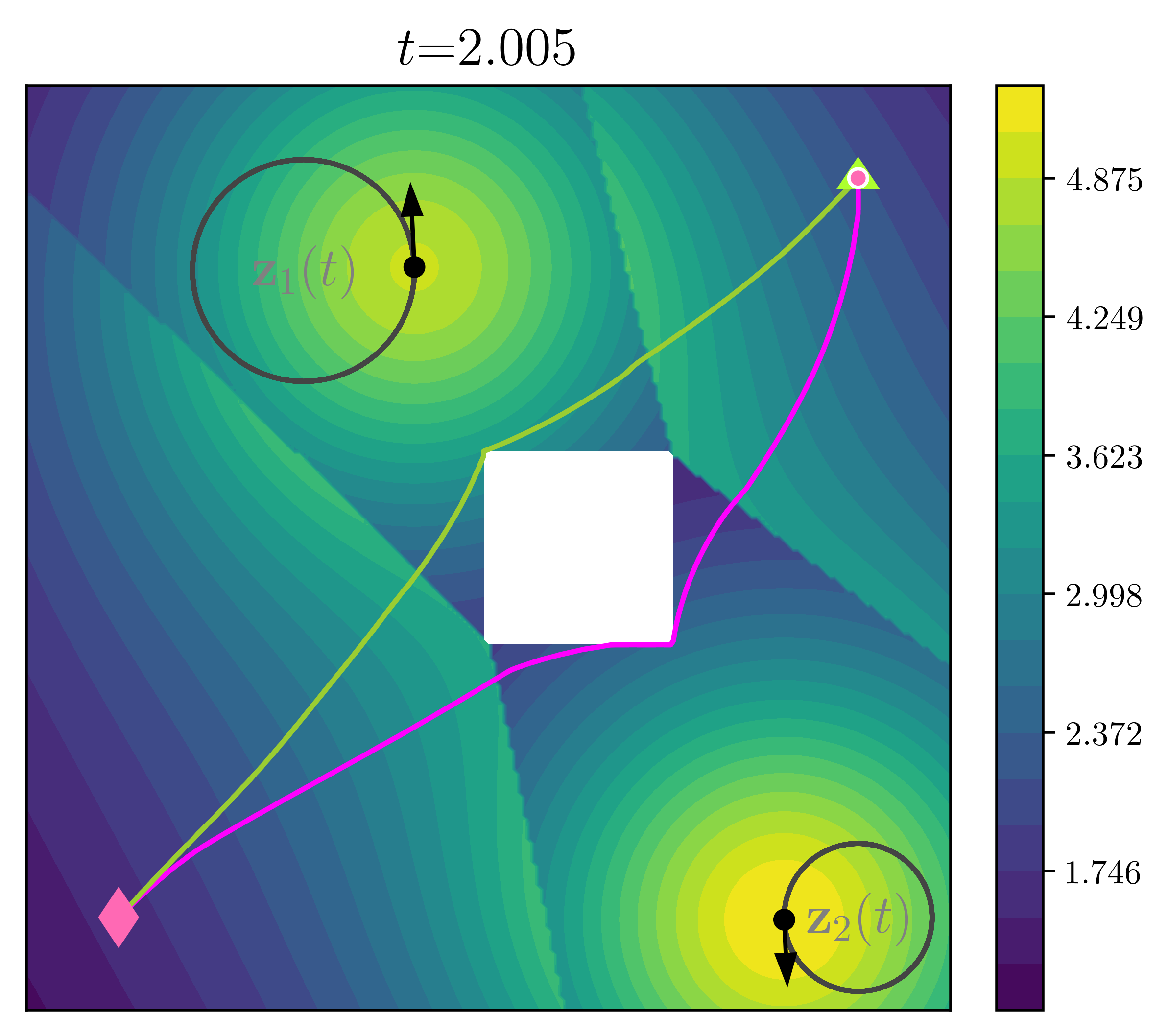}
\end{array}
$
\caption{Four snapshots of the Nash Equilibrium solution for Example 2.  Observer's optimal policy is $\lambda_*=(0.48, 0.52)$, and Evader's optimal policy is $\theta_*=(0.691, 0.309),$ with two different $\lambda_*$-optimal trajectories shown, one in yellow and one in blue corresponding to the yellow and blue points in Fig. \ref{fig:PF}B. \label{fig:ex2}}
\end{figure*}

\begin{figure*}[p!]
\centering
$
\arraycolsep=1pt\def\arraystretch{0.1}
\begin{array}{cccc}
\includegraphics[height=0.18\textheight]{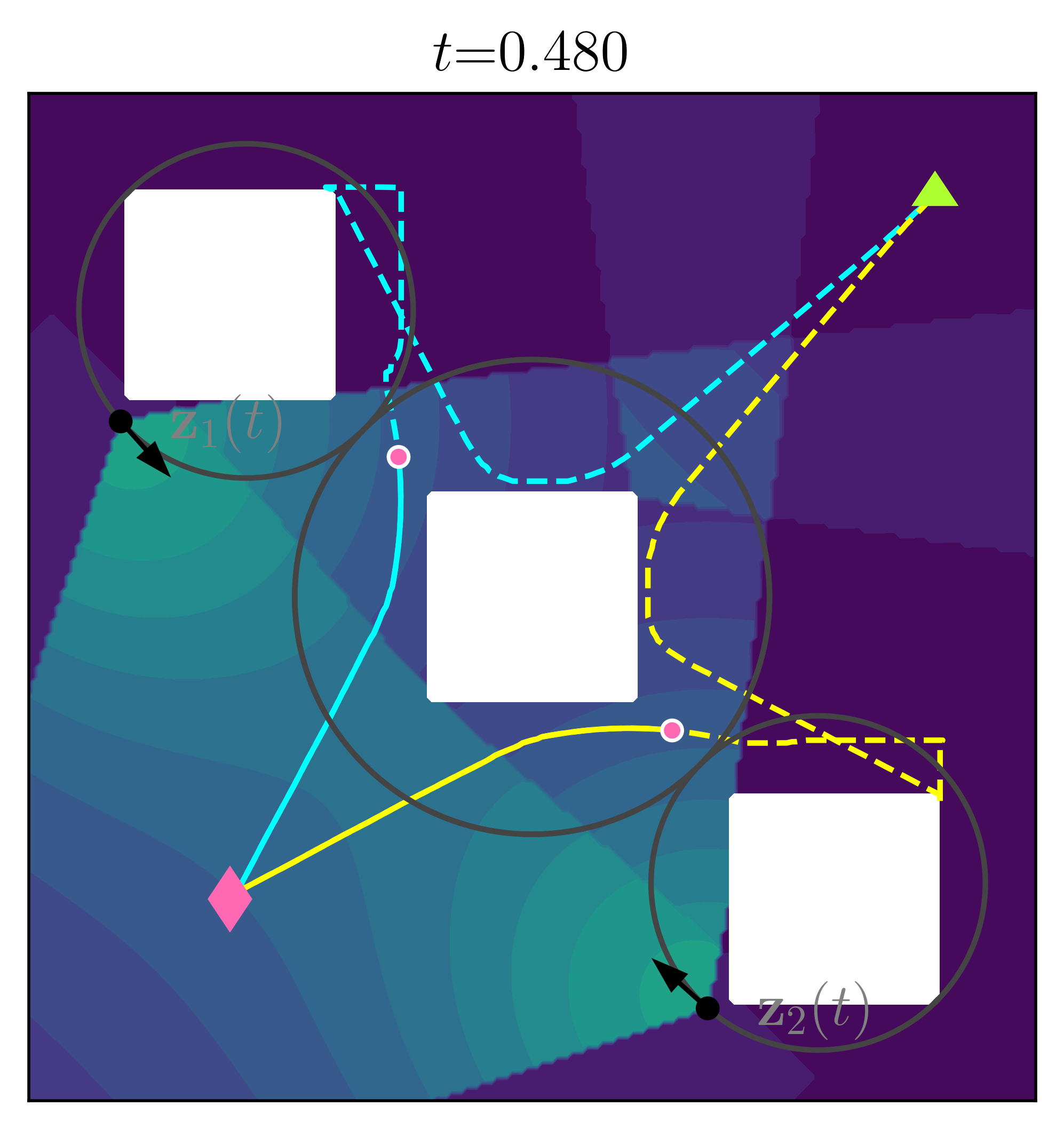} &
\includegraphics[height=0.18\textheight]{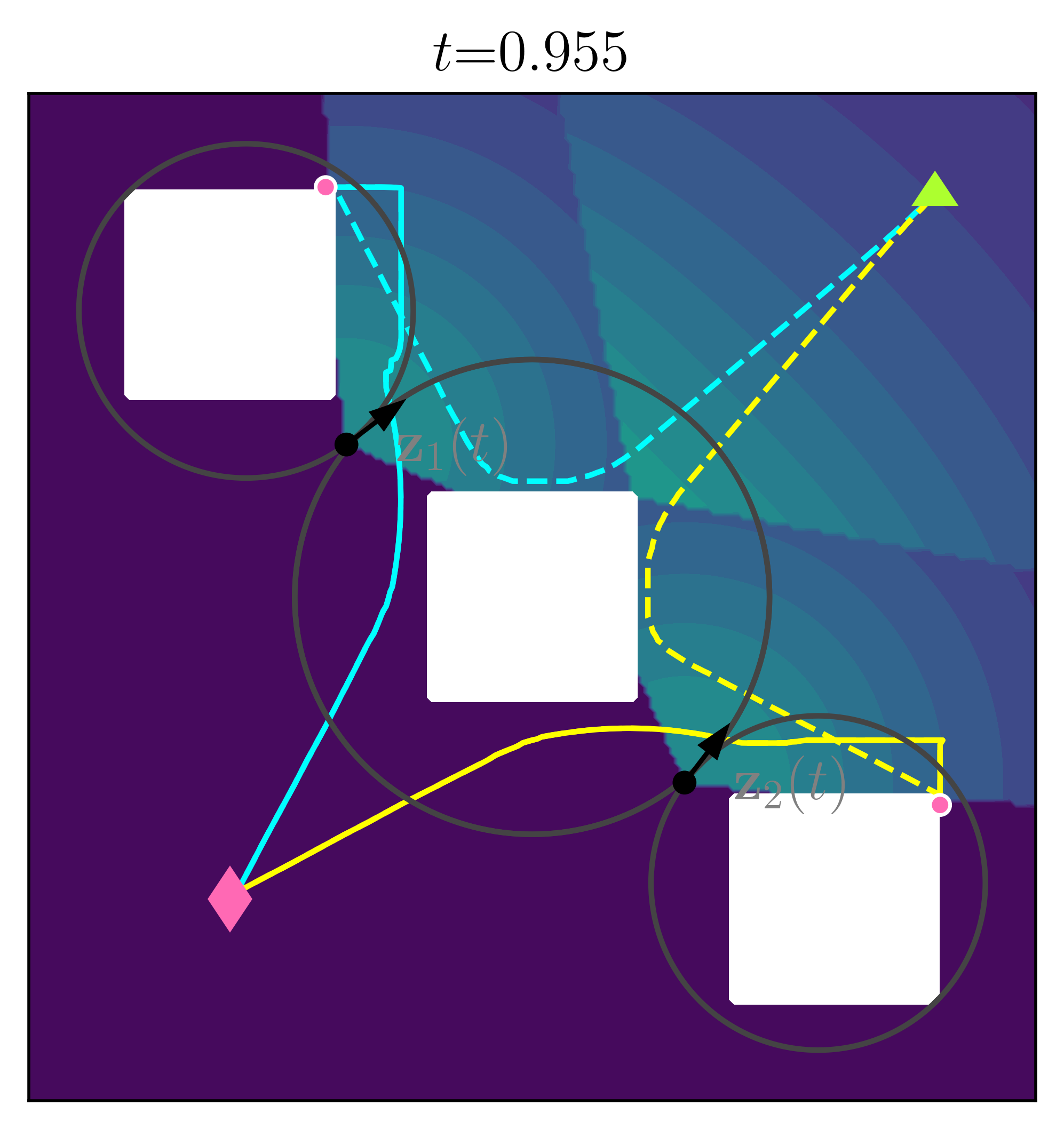} &
\includegraphics[height=0.18\textheight]{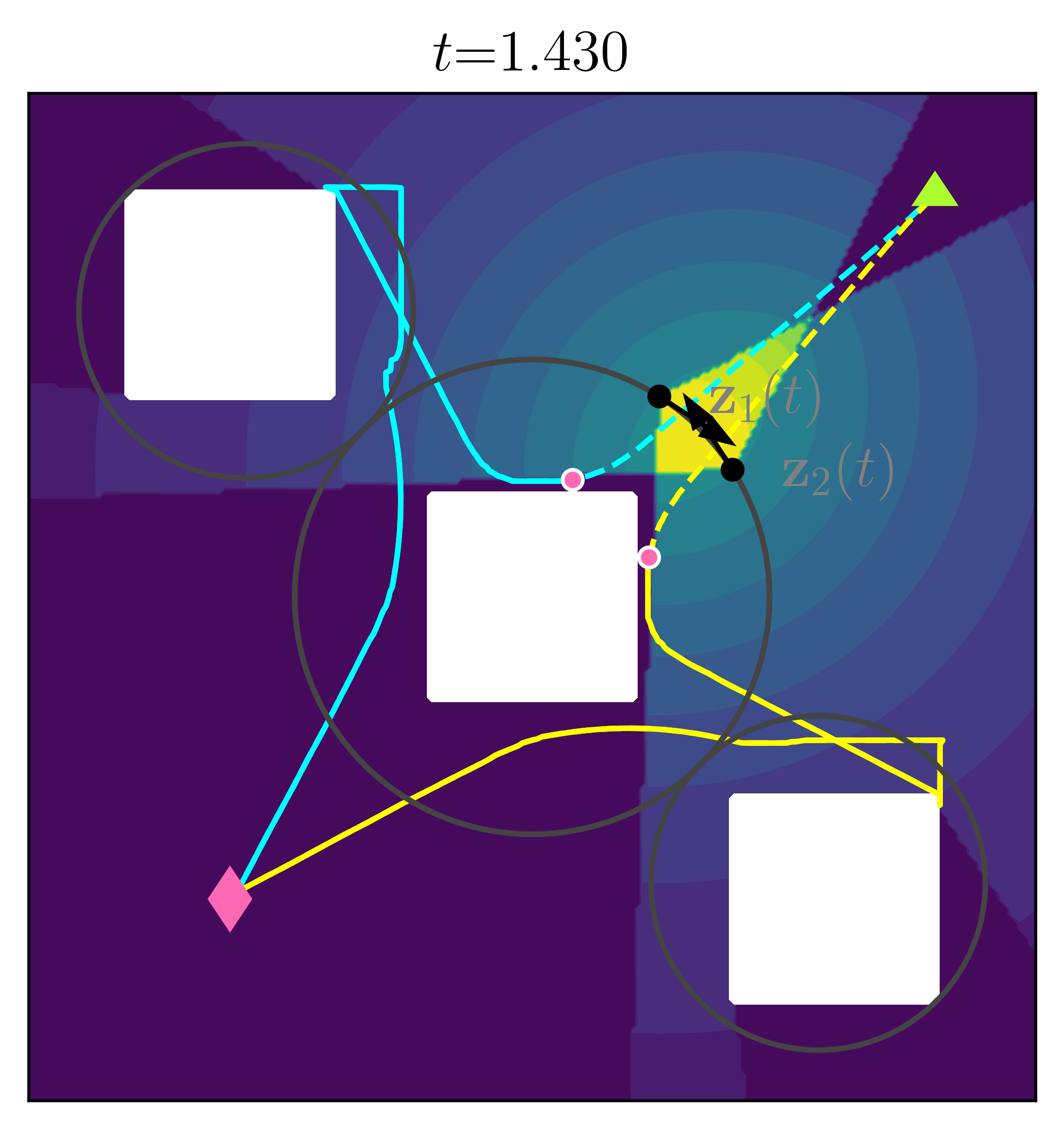} &
\includegraphics[height=0.18\textheight]{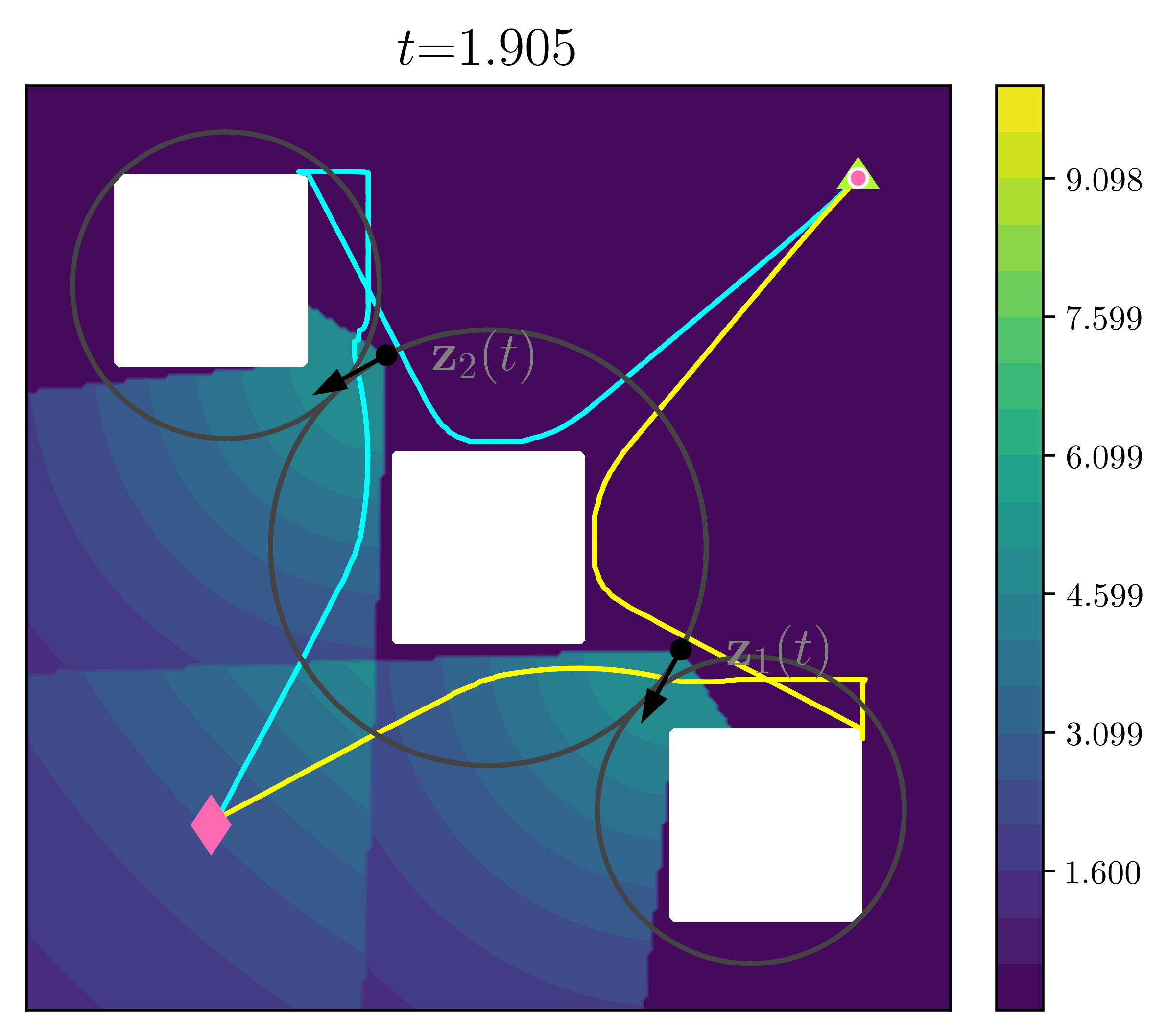}
\end{array}
$
\caption{Four snapshots of the Nash Equilibrium solution for Example 3.  Observer's optimal policy is $\lambda_*=(0.5,0.5)$, and Evader's optimal policy is $\theta_*=(0.5,0.5),$ with two different $\lambda_*$-optimal trajectories shown in cyan and yellow. \label{fig:ex3}}
\end{figure*}

\begin{figure*}[p!]
\centering
$
\arraycolsep=1pt\def\arraystretch{0.1}
\begin{array}{cccc}
\includegraphics[height=0.18\textheight]{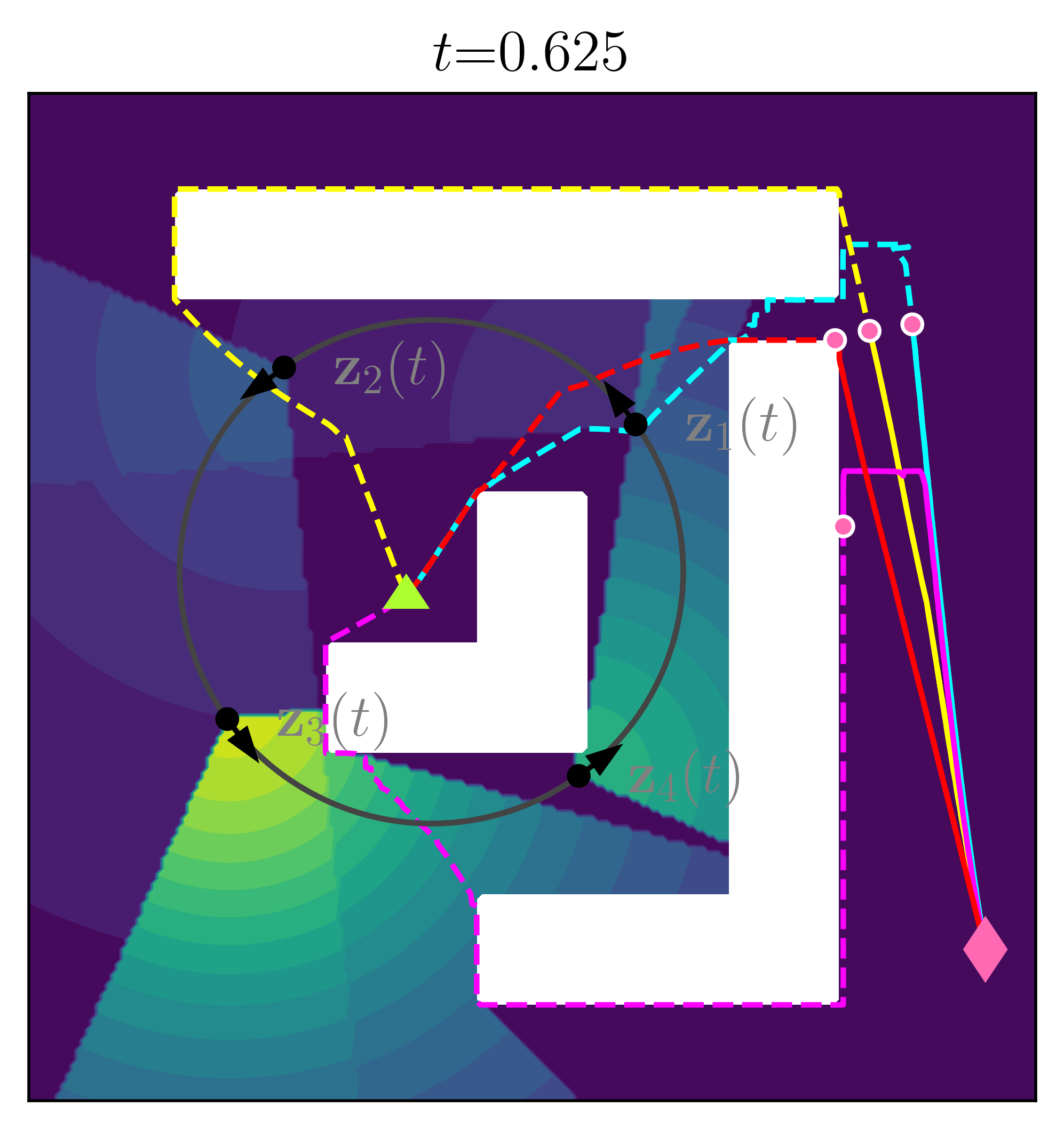} &
\includegraphics[height=0.18\textheight]{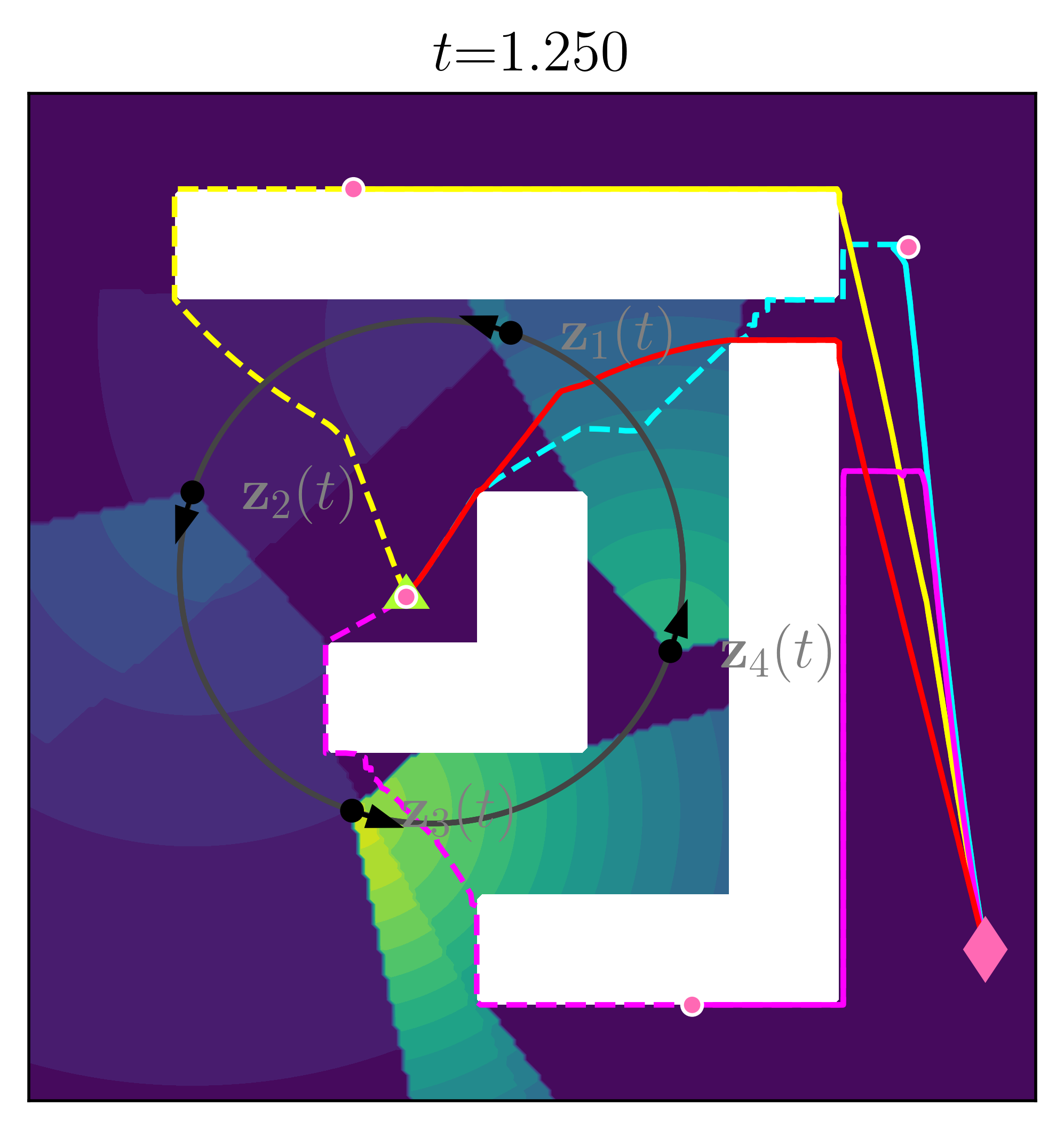} &
\includegraphics[height=0.18\textheight]{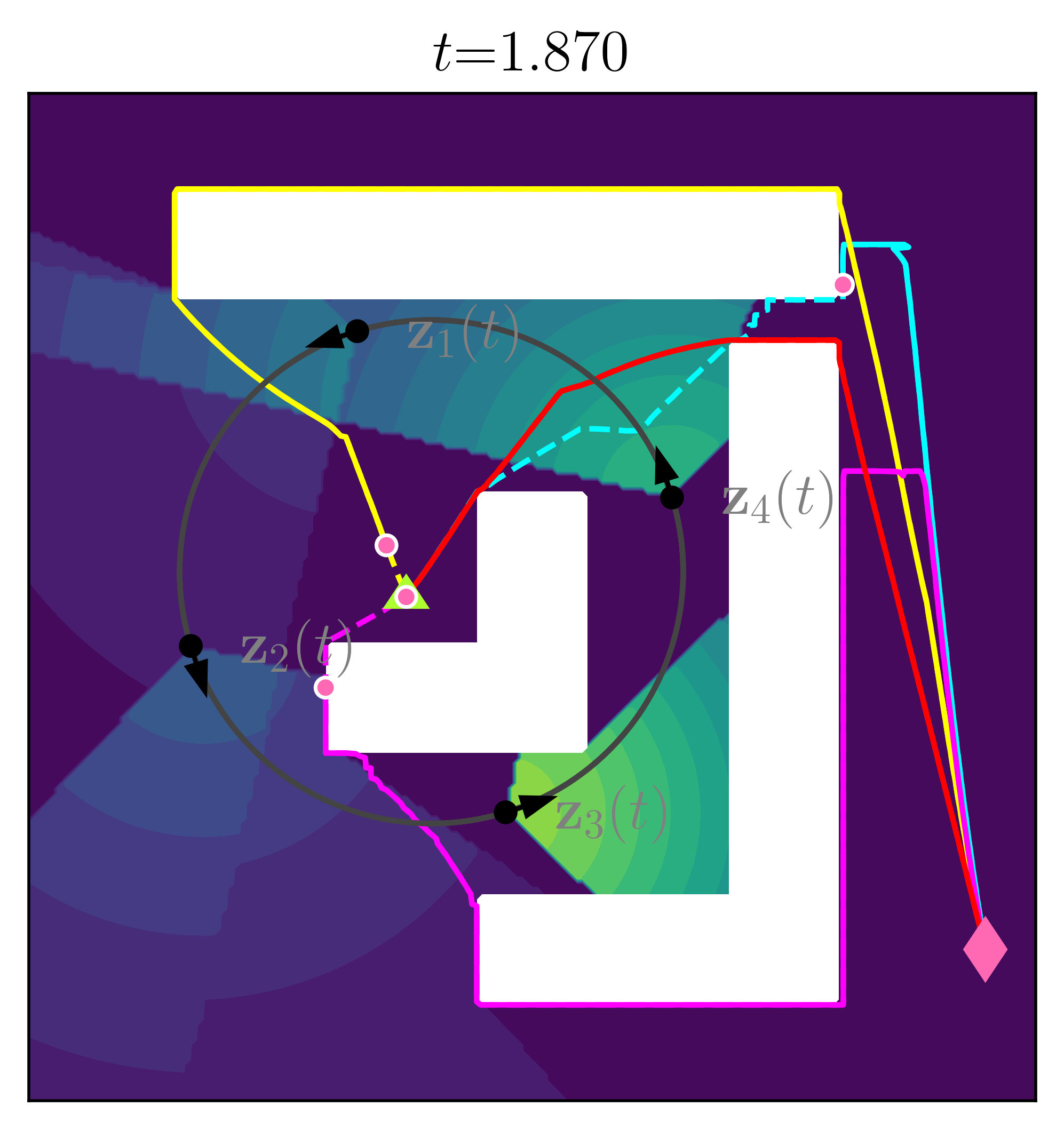} &
\includegraphics[height=0.18\textheight]{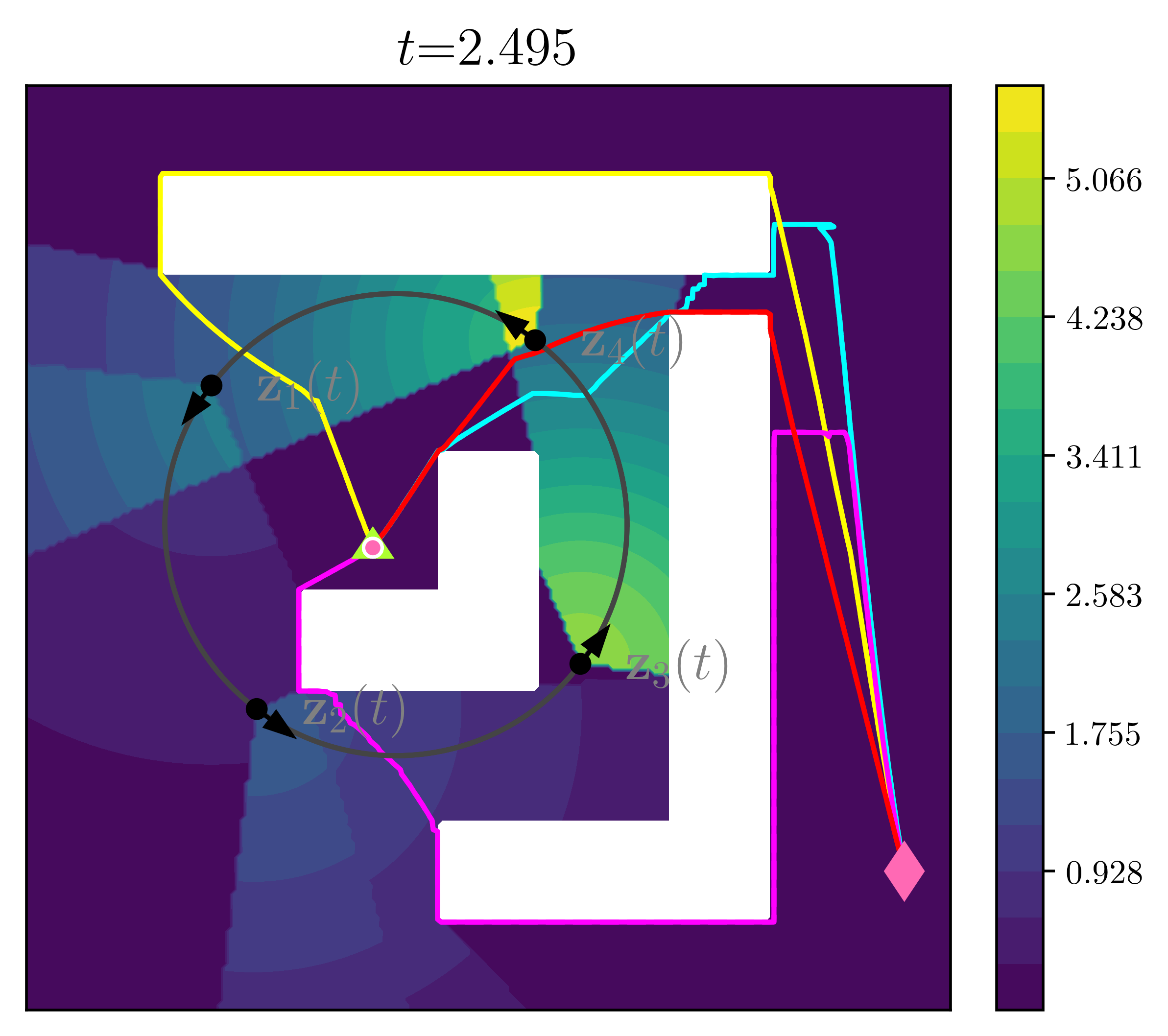}
\end{array}
$
\caption{Four snapshots of the Nash Equilibrium solution for Example 4.  Observer's optimal policy is $\lambda_*=(0.077, 0.127, 0.452, 0.344)$, and Evader's optimal policy is $\theta_*=(0.592, 0.084, 0.145, 0.180),$ with four different $\lambda_*$-optimal trajectories shown in cyan, yellow, magenta, and red. \label{fig:ex4}}
\end{figure*}

In all of our examples\footnote{
Movies 
for all of these examples can be found at \href{https://eikonal-equation.github.io/TimeDependent_SEG}{\url{https://eikonal-equation.github.io/TimeDependent_SEG}}.}, 
we assume that the domain is a unit square, sometimes containing impenetrable and occluding obstacles.
This square is discretized on a $201 \times 201$ grid ($h = 0.005$).
E's speed is uniform ($f(\bx) = 1$) and the timestep is chosen according to the CFL stability condition: $\Delta t = h.$
Our examples use observability functions $K_i(\bx,t)$ of the form found in \eqref{eq:obs}, with $K_0 = 1$
and $\sigma = 0.1$.  
In all examples, the Evader's deadline for reaching the target is $T=4,$ but all Nash equilibrium controls lead to a much earlier arrival.
(In our setting, the cumulative observability is computed up to the arrival time $T_{\ba}$ only.)

{\bf Example 1:} No obstacles, two patrol trajectories, and omnidirectional sensors (i.e., $\alpha = 2\pi$).  
Observer uses the same angular speed and moves counterclockwise on both patrol trajectories, which were already introduced in Fig.~\ref{fig:pure}.
Nash equilibrium policies are shown in Fig.~\ref{fig:ex1}.  
See also the Pareto Front in Fig.~\ref{fig:PF}(A). 

{\bf Example 2:} One obstacle, with the same sensors and patrol trajectories as Example 1.
Nash equilibrium policies are shown in Fig.~\ref{fig:ex2}. Evader's optimal policy is probabilistic, relying on two $\lambda_*$-optimal trajectories that pass above and below the obstacle.  See also the Pareto Front in Fig.~\ref{fig:PF}(B).

{\bf Example 3:} Three obstacles, two patrol trajectories with direction-restricted sensors ($\alpha = 2\pi/3$).
The patrol trajectories each follow a figure-eight pattern.
Nash equilibrium policies are shown in Fig.~\ref{fig:ex3}. 

{\bf Example 4:} Many obstacles, with 4 patrol trajectories on the same circle, but at different phases.
Observer's angular speed is $\omega = 1/2\pi$, and with direction-restricted sensors (i.e., $\alpha = 2\pi/3$). 
Nash equilibrium policies are shown in Fig.~\ref{fig:ex4}. 


\section{Conclusions}
\label{s:conclusions}
We have developed and tested a numerical method for SE games under uncertainty with isotropic Evader dynamics, a time-dependent Observer, and a direction-limited sensor.  The efficiency of our approach stems from combining existing techniques from game theory, convex optimization, and dynamic programming in continuous state/time.  The latter requires efficient numerical methods for time-dependent Eikonal PDEs.  Incorporating more realistic models of Evader dynamics will require solving more general Hamilton-Jacobi PDEs.
While our current implementation relies on an explicit first-order upwind scheme, the use of higher-order accurate discretizations \cite{Falcone2013,shu2007high} would be desirable in the future.
Additional speed improvements might be also attained through a time-dependent version of causal domain restriction \cite{CausalDomRestrict}.  All numerical examples presented above assumed that the Observer is a robotic platform with the direction of sensor/camera fixed relative to the Observer's heading.  However, the same methods would work almost without change for panning/tilting cameras, Observers with multiple sensors, and more sophisticated observability models in the Observer's field of view \cite{wang2013sparse}.  It is also relatively easy to extend the approach to a centralized planner controlling several Evaders, but if all of them have separate targets, the cost of each iteration will increase accordingly.  For a large group of Evaders, one will need to rely on distributed optimal control techniques or on Mean Field Game models, if these Evaders engage in selfish/independent trajectory planning.

\noindent
{\bf Acknowledgements:} The authors are indebted to Marc Gilles, whose algorithms and implementation for the stationary case served as a foundation for this work.  The authors are also grateful to Casey Garner and Tristan Reynoso for their valuable help in implementation and testing during the summer REU-2018 program at Cornell University.

\addtolength{\textheight}{0cm}




\bibliographystyle{IEEEtran}
\bibliography{IEEEabrv,mybibfile}

\end{document}